\documentclass[fleqn,11pt,a4paper,leqno]{article}

\usepackage{graphicx}
\usepackage{graphics}
\usepackage{latexsym}
\usepackage{amssymb,amsmath}
\usepackage{multicol}
\usepackage{bussproofs}
\usepackage{enumerate}
\usepackage{xcolor}
\usepackage[latin1]{inputenc}

\numberwithin{equation}{section}

\title{Proofs of some Propositions of the semi-Intuitionistic Logic with Strong Negation }

\author{\sc Juan Manuel Cornejo \and \sc Ignacio Viglizzo}

\date{}

\newtheorem{definicion}{Definition}[section]

\newtheorem{theorem}[definicion]{Theorem}

\newtheorem{lema}[definicion]{Lemma}

\newenvironment{Proof}{\noindent\bf Proof \rm}{$\hfill \square$}

\newcounter{numeroaxioma}

\newcommand{\il}{\to} 
\newcommand{\inl}{\to_N} 
\newcommand{\brig}{\Rightarrow}

\newcommand{\N}{\mathcal{N}} 
\newcommand{\SNlog}{\mathcal{SN}} 

\newcommand{\cl}{\vdash} 

\newcommand{\variables}{Var} 
\newcommand{\formulas}[1]{Fm_{\mathbf{#1}}} 
\newcommand{\algformulas}[1]{\sf{Fm}_{\mathbf{#1}}} 

\newcommand{\debaj}[2]{ #1 \brig #2}

\begin{document}

\maketitle

\begin{abstract}
   We offer the proofs that complete our article introducing the  propositional calculus called  semi-intuitionistic logic with strong negation.
	\end{abstract}

\section{Introduction}

The Proofs of Lemmas 3.1 and 3.3 (Lemmas \ref{propiedades_Calculo} and \ref{propiedades_Calculo_conDeduccion_para_reticulado} below) were left out of our article \cite{cornejo17semiIntutionistic}. We detail them here, together with the necessary axioms for the semi-Intuitionistic logic with strong negation. Please refer to that article for motivation and more results on this calculus.

\section{Semi-intuitionistic logic with strong negation}

 A {\it logical language} $\mathbf L$, as defined in \cite{Font03survey}, is  a set of connectives, each with a fixed arity $n \geq 0$. For a countably infinite set $\variables$ of propositional variables, the {\it formulas} of the logical language $\mathbf L$
 are inductively defined as usual.  

A {\it logic}, in the language $\mathbf L$, is a pair $\mathcal L = \langle \formulas{L}, \cl_\mathcal{L} \rangle$ where $\formulas{L}$ is the set of formulas and $\cl_\mathcal{L}$ is a substitution-invariant consequence relation on $\formulas{L}$.  As usual, the set $\formulas{L}$ may also be endowed with an algebraic structure, just by regarding the connectives of the language as operation symbols. The resulting algebra  is the {\it algebra of formulas}, denoted by $\algformulas{L}$. The finitary logic is presented by means of their ``Hilbert style'' sets of axioms and inferences rules.
 
 \medskip
 
We define {\it semi Intuitionistic logic with strong negation}  $\SNlog$ over the language $\mathbf L=\{\top, \sim, \wedge, \vee, \il\}$ in terms of the following set of axiom schemata, in which we use the following definitions: 
\begin{itemize}
	\item $\alpha \inl \beta : = \alpha \il (\alpha \wedge \beta)$,
	\item $\alpha \brig \beta := (\alpha \inl \beta) \wedge (\sim\beta \inl \sim \alpha)$.
\end{itemize}


\medskip

\begin{enumerate}[$({A}1)$]
	
	\item $(\alpha \inl \beta) \inl ((\beta \inl \gamma) \inl (\alpha \inl \gamma))$, \label{axioma_transitividadNelson} 
	\item $(\alpha \inl \beta) \inl ((\alpha \inl \gamma) \inl (\alpha \inl (\beta \wedge \gamma)))$, \label{axioma_mayor_cota_inferior} 
	\item $(\alpha \wedge \beta) \inl \alpha$, \label{axioma_infimo_izquierda}
	\item $(\alpha \wedge \beta) \inl \beta$, \label{axioma_infimo_derecha}
	\item $\alpha \inl (\alpha \vee \beta)$, \label{axioma_supremo_izquierda} %
	\item $\beta \inl (\alpha \vee \beta)$,  \label{axioma_supremo_derecha} %
	\item $\sim (\alpha \vee \beta) \inl \sim\alpha$, \label{axioma_supremo_negado_izquierda}
	\item $\sim (\alpha \vee \beta) \inl \sim\beta$, \label{axioma_supremo_negado_derecha}
	\item $(\alpha \inl \gamma) \inl ((\beta \inl \gamma) \inl ((\alpha \vee \beta) \inl \gamma))$, \label{axioma_menor_cota_superior}
	\item \label{axioma_menor_cota_superior_negado} $(\sim\alpha \inl \sim\beta) \inl ((\sim\alpha \inl\sim \gamma) \inl (\sim\alpha \inl \sim(\beta \vee \gamma)))$, 
	\item $\alpha \brig (\sim \sim \alpha)$, \label{axioma_doble_neg_a_derecha}
	\item $(\sim \sim \alpha) \brig \alpha$, \label{axioma_doble_neg_a_izquierda}
	\item $(\alpha \inl \beta) \inl [(\beta \inl \alpha) \inl [(\alpha \il \gamma) \inl (\beta \il \gamma)]]$,  \label{axioma_buena_def_implica_der}
	\item $(\alpha \inl \beta) \inl [(\beta \inl \alpha) \inl [(\gamma \il \alpha) \inl (\gamma \il \beta)]]$, \label{axioma_buena_def_implica_izq} 
	\item $[(\alpha \wedge \beta) \inl \gamma] \brig [\alpha \inl (\beta \inl \gamma)]$,  \label{axioma_InfimoAImplicacion}
	

	\item $\debaj{(\sim (\alpha \wedge \beta))}{(\sim \alpha \vee \sim \beta)}$, \label{axioma_distribuye_neg_infimo1}
	\item $\debaj{(\sim \alpha \vee \sim \beta)}{(\sim (\alpha \wedge \beta))}$, \label{axioma_distribuye_neg_infimo2}
	\item \label{axioma_implica_infimo_dosVariables2} $\debaj{(\alpha \wedge (\sim \alpha \vee \beta))}{(\alpha \wedge (\alpha \inl \beta))}$, 
	\item $\debaj{(\alpha \inl (\beta \inl \gamma))}{((\alpha \wedge \beta) \inl \gamma)}$,  \label{axioma_InfimoAImplicacionVuelta}
	\item $(\sim (\alpha \il \beta)) \inl (\alpha \wedge \sim \beta)$,  \label{axioma_paraSN1}
	\item $(\alpha \wedge \sim \beta) \inl (\sim (\alpha \il \beta))$,  \label{axioma_paraSN2}
	\item $[\sim (\alpha \wedge ((\gamma \wedge \alpha) \vee (\beta \wedge \alpha)))] \inl [\sim (\alpha \wedge (\beta \vee \gamma))]$, \label{axioma_para_reticulado_neg}
	\item \label{axioma_top} $\top$. 
\end{enumerate}

The only inference rule is Modus Ponens  for the implication $\inl$, which we denominate  {\it $\N$-Modus Ponens } ($\N$-MP):  $\Gamma
	\cl_{\SNlog} \phi$ and $\Gamma \cl_{\SNlog} \phi \inl \gamma$ yield
	$\Gamma \cl_{\SNlog} \gamma$.

\begin{lema} \label{propiedades_Calculo}
	Let $\Gamma \cup \{\alpha, \beta\} \subseteq \formulas{L}$. In $\SNlog$ the following properties hold:
	\begin{enumerate}[{\rm (a)}]
		\item If $\Gamma \cl \alpha$ then $\Gamma \cl \beta \inl \alpha$, \label{implicacionATeorema}
		\item $ \Gamma \cl \alpha \inl \alpha$, \label{propiedad_XimplicaX}
		\item If $ \Gamma \cl \alpha \brig \beta$ then $ \Gamma \cl \alpha \inl \beta$ and $ \Gamma \cl \sim\beta \inl \sim\alpha$, \label{brigAFactores}
		\item $\Gamma \cl \sim\alpha \inl \sim(\alpha \wedge \beta)$, \label{110716_02}  
		
			\item $\Gamma \cl \sim\beta \inl\sim (\alpha \wedge \beta)$, \label{110716_03}  
				
		\item $ \Gamma, \alpha, \alpha \brig \beta \cl \beta$,  \label{condicionIL4}
		\item If $ \Gamma \cl \alpha \brig \beta$ and  $ \Gamma \cl \alpha$ then $ \Gamma \cl \beta$, \label{MPsobreBrigN}
		\item If $\Gamma \cl \alpha$ and $\Gamma \cl \beta$ then $\Gamma \cl \alpha \wedge \beta$, \label{021115_01}
		\item $\Gamma \cl \alpha \wedge \beta \brig \alpha$ and $\Gamma \cl \alpha \wedge \beta \brig \beta$, \label{260216_37}
		\item $\Gamma \cl \alpha \brig \alpha \vee \beta$ and $\Gamma \cl \alpha  \brig \beta \vee \alpha$, \label{120416_01}
		\item $\Gamma \cl \alpha \brig \alpha$, \label{condicionIL1}
		\item If $\Gamma \cl \alpha \brig \beta$ and $\Gamma \cl \beta \brig \gamma$ then $\Gamma \cl \alpha \brig \gamma$, \label{propiedad_Transitividad}
		\item $\Gamma, \alpha \brig \beta, \beta \brig \gamma \cl \alpha \brig \gamma$, \label{condicionIL2}
				\item $\Gamma \cl \alpha \inl \beta$ then $\Gamma \cl (\gamma \wedge \alpha) \inl (\gamma \wedge \beta)$ and $\Gamma \cl (\alpha \wedge \gamma) \inl (\beta \wedge \gamma)$, \label{280416_01}
				
				\item $\Gamma \cl \alpha \inl \beta$ then $\Gamma \cl (\gamma \vee \alpha) \inl (\gamma \vee \beta)$ and $\Gamma \cl (\alpha \vee \gamma) \inl (\beta \vee \gamma)$,  \label{280416_14}
				
		\item $\Gamma \cl (\alpha \vee \beta) \inl (\beta \vee \alpha)$, \label{080716_42} 
		\item $\Gamma \cl (\alpha \land \beta) \inl (\beta \land \alpha)$, \label{290916_00} 
		\item $\Gamma, \alpha \brig \beta \cl (\alpha \vee \gamma) \brig (\beta \vee \gamma)$, \label{111215_42}
		\item $\Gamma, \alpha \brig \beta \cl (\gamma \vee \alpha) \brig (\gamma \vee \beta)$,  \label{111215_43}
		\item $\Gamma, \alpha \brig \beta, \gamma \brig t \cl (\alpha \vee \gamma) \brig (\beta \vee t)$, \label{condicionIL3b}
		\item $\Gamma,  \beta \brig \alpha \cl (\sim \alpha) \brig (\sim \beta)$, \label{condicionIL3c}
		\item $\Gamma \cl (\sim (\alpha \il \beta)) \inl (\sim (\alpha \inl \beta))$. \label{151116_11} 
	\end{enumerate}
\end{lema}


\begin{Proof} 
	\begin{itemize}
		\item[(\ref{implicacionATeorema})]
		\begin{enumerate}[1.]
			\item $ \Gamma \cl [(\alpha \wedge \beta) \inl \alpha] \brig [\alpha \inl (\beta \inl \alpha)]$ by axiom $(A{\ref{axioma_InfimoAImplicacion}})$. \label{250216_01}
			\item $ \Gamma \cl [[(\alpha \wedge \beta) \inl \alpha] \brig [\alpha \inl (\beta \inl \alpha)]] \inl [[(\alpha \wedge \beta) \inl \alpha] \inl [\alpha \inl (\beta \inl \alpha)]]$  by axiom $(A{\ref{axioma_infimo_izquierda}})$. \label{250216_02}
			\item $ \Gamma \cl [[(\alpha \wedge \beta) \inl \alpha] \inl [\alpha \inl (\beta \inl \alpha)]$  by ($\N$-MP) applied to  \ref{250216_01} and \ref{250216_02}. \label{250216_03}
			\item $ \Gamma \cl (\alpha \wedge \beta) \inl \alpha$ by axiom $(A{\ref{axioma_infimo_izquierda}})$.  \label{250216_04}
			\item $ \Gamma \cl \alpha \inl (\beta \inl \alpha)$ by ($\N$-MP) applied to \ref{250216_03} and \ref{250216_04}. \label{250216_05}
			\item $ \Gamma \cl \alpha$ by hypothesis. \label{250216_06}
			\item $ \Gamma \cl \beta \inl \alpha$  by ($\N$-MP) applied to \ref{250216_05} and \ref{250216_06}.
		\end{enumerate}

		\item[(\ref{propiedad_XimplicaX})] Let $\phi$ be any axiom of $\SNlog$.
		\begin{enumerate}[1.]
			\item $\Gamma \cl \phi$. \label{211215_08}
			\item $\Gamma \cl [(\phi \wedge \alpha) \inl \alpha] \brig [\phi \inl (\alpha \inl \alpha)]$
			by axiom  $(A{\ref{axioma_InfimoAImplicacion}})$. \label{211215_04}
			\item $\Gamma \cl \{[(\phi \wedge \alpha) \inl \alpha] \brig [\phi \inl (\alpha \inl \alpha)]\} \inl \{[(\phi \wedge \alpha) \inl \alpha] \inl [\phi \inl (\alpha \inl \alpha)]\}$ by axiom $(A{\ref{axioma_infimo_izquierda}})$. \label{211215_05}
			\item $\Gamma \cl [(\phi \wedge \alpha) \inl \alpha] \inl [\phi \inl (\alpha \inl \alpha)]$  by ($\N$-MP) applied to \ref{211215_04} and \ref{211215_05}. \label{211215_06}
			\item $\Gamma \cl (\phi \wedge \alpha) \inl \alpha$  by axiom $(A{\ref{axioma_infimo_derecha}})$. \label{211215_07}
			\item $\Gamma \cl \phi \inl (\alpha \inl \alpha)$  by ($\N$-MP) applied to \ref{211215_06} and \ref{211215_07}. \label{211215_09}
			\item $\Gamma \cl \alpha \inl \alpha$  by ($\N$-MP) applied to \ref{211215_08} and \ref{211215_09}.
		\end{enumerate}

		\item[(\ref{brigAFactores})] 
		\begin{enumerate}[1.]
			\item $ \Gamma \cl \alpha \brig \beta$ by hypothesis.
			\item $ \Gamma \cl (\alpha \inl \beta) \wedge (\sim\beta \inl \sim\alpha)$ by the definition of $\brig$.
			\item $ \Gamma \cl [(\alpha \inl \beta) \wedge (\sim\beta \inl \sim\alpha)] \inl (\alpha \inl \beta)$ by axiom $(A{\ref{axioma_infimo_izquierda}})$.
			\item $ \Gamma \cl \alpha \inl \beta$ by  ($\N$-MP).
			\item $ \Gamma \cl [(\alpha \inl \beta) \wedge (\sim\beta \inl \sim\alpha)] \inl (\sim\beta \inl\sim \alpha)$ by axiom $(A{\ref{axioma_infimo_derecha}})$.
			\item $ \Gamma \cl \sim\beta \inl \sim\alpha$ by  ($\N$-MP).
		\end{enumerate}

		\item[ (\ref{110716_02})] 
		\begin{enumerate}
			\item 	$\Gamma \cl (\sim \alpha) \inl (\sim \alpha \vee \sim \beta)$ by axiom $(A{\ref{axioma_supremo_izquierda}})$.  \label{190416_13}
			\item 	$\Gamma \cl (\sim \alpha \vee \sim \beta) \brig (\sim (\alpha \wedge \beta))$ by axiom $(A{\ref{axioma_distribuye_neg_infimo2}})$. \label{270417_01}
			\item 	$\Gamma \cl (\sim \alpha \vee \sim \beta) \inl (\sim (\alpha \wedge \beta))$ by part (\ref{brigAFactores}) applied to \ref{270417_01}. \label{190416_14}
			\item 	$\Gamma \cl (\sim \alpha) \inl (\sim (\alpha \wedge \beta))$  by axiom $(A{\ref{axioma_transitividadNelson}})$ and ($\N$-MP) applied to \ref{190416_13} and \ref{190416_14}.
		\end{enumerate}

		\item[(\ref{110716_03})] 
		\begin{enumerate}
			\item 	$\Gamma \cl (\sim \beta) \inl (\sim \alpha \vee \sim \beta)$ by axiom $(A{\ref{axioma_supremo_derecha}})$.  \label{190416_15}
			\item 	$\Gamma \cl (\sim \alpha \vee \sim \beta) \brig (\sim (\alpha \wedge \beta))$ by axiom $(A{\ref{axioma_distribuye_neg_infimo2}})$. \label{270417_02}
			\item 	$\Gamma \cl (\sim \alpha \vee \sim \beta) \inl (\sim (\alpha \wedge \beta))$ by part (\ref{brigAFactores}) applied to \ref{270417_02}. \label{190416_16}
			\item 	$\Gamma \cl (\sim \beta) \inl (\sim (\alpha \wedge \beta))$  by axiom $(A{\ref{axioma_transitividadNelson}})$ and ($\N$-MP) applied to \ref{190416_13} and \ref{190416_14}.
		\end{enumerate}
		
		\item[(\ref{condicionIL4})] 
		\begin{enumerate}[1.]
			\item $ \Gamma, \alpha, \alpha \brig \beta \cl \alpha \brig \beta$. 
			\item $ \Gamma, \alpha, \alpha \brig \beta \cl \alpha \inl \beta$  by part (\ref{brigAFactores}). \label{111215_68}
			\item $ \Gamma, \alpha, \alpha \brig \beta \cl \alpha$  \label{111215_69}
			\item $ \Gamma, \alpha, \alpha \brig \beta \cl \beta$  by ($\N$-MP) applied to \ref{111215_68} and \ref{111215_69}.
		\end{enumerate}
		
		\item[(\ref{MPsobreBrigN})] 
		\begin{enumerate}[1.]
			\item $ \Gamma \cl \alpha \brig \beta$ by hypothesis.
			\item $ \Gamma \cl \alpha \inl \beta$ by part (\ref{brigAFactores}). \label{270417_03}
			\item $ \Gamma \cl \alpha$ by hypothesis. \label{270417_04}
			\item $ \Gamma \cl \beta$ by  ($\N$-MP) applied to \ref{270417_03} and \ref{270417_04}.
		\end{enumerate}
		
		\item[(\ref{021115_01})] 
		\begin{enumerate}[1.]
			\item $\Gamma \cl [(\beta \wedge \alpha) \inl \beta] \brig [\beta \inl (\alpha \inl \beta)]$  by axiom $(A{\ref{axioma_InfimoAImplicacion}})$. \label{301015_01}
			\item $\Gamma \cl (\beta \wedge \alpha) \inl \beta$ by axiom $(A{\ref{axioma_infimo_izquierda}})$. \label{301015_02}
			\item $\Gamma \cl \beta \inl (\alpha \inl \beta)$ by (\ref{MPsobreBrigN}) applied to \ref{301015_01} and \ref{301015_02}. \label{301015_03}
			\item $\Gamma \cl \beta$ by hypothesis. \label{301015_04}
			\item $\Gamma \cl \alpha \inl \beta$ by ($\N$-MP) applied to \ref{301015_03} and \ref{301015_04}. \label{301015_07}
			\item $\Gamma \cl (\alpha \inl \alpha) \inl ((\alpha \inl \beta) \inl (\alpha \inl (\alpha \wedge \beta)))$ by axiom $(A{\ref{axioma_mayor_cota_inferior}})$. \label{301015_05}
			\item $\Gamma \cl  \alpha \inl \alpha$  by part (\ref{propiedad_XimplicaX}). \label{301015_06}
			\item  $\Gamma \cl  (\alpha \inl \beta) \inl (\alpha \inl (\alpha \wedge \beta))$  by ($\N$-MP) applied to \ref{301015_05} and \ref{301015_06}. \label{301015_08}
			\item  $\Gamma \cl  \alpha \inl (\alpha \wedge \beta)$  by ($\N$-MP) applied to \ref{301015_07} and \ref{301015_08}. \label{301015_09}
			\item $\Gamma \cl \alpha$ by hypothesis. \label{301015_10}
			\item $\Gamma \cl \alpha \wedge \beta$  by ($\N$-MP) applied to \ref{301015_09} and \ref{301015_10}.
		\end{enumerate}
		
		\item[(\ref{260216_37})] Follows immediately from items (\ref{021115_01}) and (\ref{110716_02}),  axioms $(A{\ref{axioma_infimo_izquierda}})$, $(A{\ref{axioma_infimo_derecha}})$
		and item (\ref{110716_03}).
		
		\item[(\ref{120416_01})] Is a direct consequence of  item (\ref{021115_01}) and axioms $(A{\ref{axioma_supremo_izquierda}})$, $(A{\ref{axioma_supremo_derecha}})$, $(A{\ref{axioma_supremo_negado_izquierda}})$, and $(A{\ref{axioma_supremo_negado_derecha}})$.
		
		\item[(\ref{condicionIL1})] 
		\begin{enumerate}[1.]
			\item $ \Gamma \cl \alpha \inl \alpha$ by part (\ref{propiedad_XimplicaX}).
			\item $ \Gamma \cl (\sim \alpha) \inl (\sim \alpha)$ by part (\ref{propiedad_XimplicaX}).
			\item $ \Gamma \cl (\alpha \inl \alpha) \wedge ((\sim \alpha) \inl (\sim \alpha))$ by part (\ref{021115_01}).
		\end{enumerate}
		
		\item[(\ref{propiedad_Transitividad})] 
		
		\begin{enumerate}[1.]
			\item $ \Gamma  \cl \alpha \brig \beta$ \label{100223_01} by hypothesis.
			\item $ \Gamma  \cl \alpha \inl \beta$   by 			\ref{100223_01} and (\ref{brigAFactores}). 			 \label{100223_05}
			\item $ \Gamma  \cl \beta \brig \gamma$ \label{100223_03} by hypothesis.
			\item $ \Gamma  \cl \beta \inl \gamma$   by  	  \ref{100223_03} and (\ref{brigAFactores}).
			\label{100223_07}
			\item $ \Gamma  \cl (\alpha \inl \beta) \inl ((\beta \inl \gamma) \inl (\alpha \inl \gamma))$  by axiom $(A{\ref{axioma_transitividadNelson}})$. \label{100223_06}
			\item $ \Gamma  \cl (\beta \inl \gamma) \inl (\alpha \inl \gamma)$  by ($\N$-MP) applied to \ref{100223_05} and \ref{100223_06}. \label{100223_08}
			\item $ \Gamma  \cl \alpha \inl \gamma$  by ($\N$-MP) applied to \ref{100223_07} and \ref{100223_08}. \label{100516_04}
			\item $ \Gamma  \cl (\alpha \brig \beta) \inl (\sim\beta \inl \sim\alpha)$ by axiom $(A{\ref{axioma_infimo_derecha}})$. \label{100223_09}
			\item $ \Gamma  \cl \sim\beta \inl \sim\alpha$   by ($\N$-MP) applied to \ref{100223_01} and \ref{100223_09}. \label{100223_13}
			\item $ \Gamma  \cl (\beta \brig \gamma) \inl (\sim\gamma \inl\sim \beta)$ by axiom $(A{\ref{axioma_infimo_derecha}})$. \label{100223_10}
			\item $ \Gamma  \cl \sim\gamma \inl \sim\beta$   by ($\N$-MP) applied to \ref{100223_03} and \ref{100223_10}. \label{100223_11}
			\item $ \Gamma  \cl (\sim\gamma \inl \sim\beta) \inl ((\sim\beta \inl\sim \alpha) \inl (\sim\gamma \inl\sim \alpha))$  by axiom $(A{\ref{axioma_transitividadNelson}})$. \label{100223_12}
			\item $ \Gamma  \cl (\sim\beta \inl \sim\alpha) \inl (\sim\gamma \inl \sim\alpha)$  by ($\N$-MP) applied to \ref{100223_11} and \ref{100223_12}. \label{100223_14}
			\item $ \Gamma  \cl \sim\gamma \inl \sim\alpha$  by ($\N$-MP) applied to \ref{100223_13} and \ref{100223_14}. \label{100516_03}
			\item $ \Gamma  \cl \alpha \brig \gamma$ by part (\ref{021115_01}) applied to \ref{100516_04} and \ref{100516_03}.
			
		\end{enumerate}

		\item[(\ref{condicionIL2})] 
		
		\begin{enumerate}[1.]
			\item $ \Gamma, \alpha \brig \beta, \beta \brig \gamma \cl \alpha \brig \beta$ 
			\item $ \Gamma, \alpha \brig \beta, \beta \brig \gamma \cl \beta \brig \gamma$ 
			\item $ \Gamma, \alpha \brig \beta, \beta \brig \gamma \cl \alpha \brig \gamma$   by part (\ref{propiedad_Transitividad}). 
		\end{enumerate}
		
		\item[(\ref{280416_01})]
		\begin{enumerate}
			\item 	$\Gamma \cl ((\gamma \wedge \alpha) \inl \gamma) \inl [((\gamma \wedge \alpha) \inl \beta) \inl [(\gamma \wedge \alpha) \inl (\gamma \wedge \beta)]]$  by axiom $(A{\ref{axioma_mayor_cota_inferior}})$. \label{280416_02}
			\item 	$\Gamma \cl (\gamma \wedge \alpha) \inl \gamma$   by axiom $(A{\ref{axioma_infimo_izquierda}})$.  \label{280416_03}
			\item 	$\Gamma \cl ((\gamma \wedge \alpha) \inl \beta) \inl [(\gamma \wedge \alpha) \inl (\gamma \wedge \beta)]$   by ($\N$-MP) applied to \ref{280416_02} and \ref{280416_03}. \label{280416_06}
			\item 	$\Gamma \cl (\gamma \wedge \alpha) \inl \alpha$   by axiom $(A{\ref{axioma_infimo_derecha}})$ \label{280416_04}
			\item 	$\Gamma \cl \alpha \inl \beta$ by hypothesis. \label{280416_05}
			\item 	$\Gamma \cl (\gamma \wedge \alpha) \inl \beta$  by axiom $(A{\ref{axioma_transitividadNelson}})$ and ($\N$-MP) applied to \ref{280416_04} and \ref{280416_05}. \label{280416_07}
			\item 	$\Gamma \cl (\gamma \wedge \alpha) \inl (\gamma \wedge \beta)$  by ($\N$-MP) applied to \ref{280416_06} and \ref{280416_07}.
			\item 	$\Gamma \cl ((\alpha \wedge \gamma) \inl \beta) \inl [((\alpha \wedge \gamma) \inl \gamma) \inl [(\alpha \wedge \gamma) \inl (\beta \wedge \gamma)]]$  by axiom $(A{\ref{axioma_mayor_cota_inferior}})$. \label{280416_10}
			\item 	$\Gamma \cl (\alpha \wedge \gamma) \inl \alpha$   by axiom $(A{\ref{axioma_infimo_izquierda}})$. \label{280416_08}
			\item 	$\Gamma \cl (\alpha \wedge \gamma) \inl \beta$  by axiom $(A{\ref{axioma_transitividadNelson}})$ and ($\N$-MP) applied to \ref{280416_08} and \ref{280416_05}. \label{280416_11}
			\item 	$\Gamma \cl ((\alpha \wedge \gamma) \inl \gamma) \inl [(\alpha \wedge \gamma) \inl (\beta \wedge \gamma)]$  by ($\N$-MP) applied to \ref{280416_10} and \ref{280416_11}. \label{280416_12}
			\item 	$\Gamma \cl (\alpha \wedge \gamma) \inl \gamma$  by axiom $(A{\ref{axioma_infimo_derecha}})$. \label{280416_13}
			\item 	$\Gamma \cl (\alpha \wedge \gamma) \inl (\beta \wedge \gamma)$  by ($\N$-MP) applied to \ref{280416_12} and \ref{280416_13}.
		\end{enumerate}

		\item[(\ref{280416_14})]
		\begin{enumerate}
			\item 	$\Gamma \cl (\gamma \inl (\gamma \vee \beta)) \inl [(\alpha \inl (\gamma \vee \beta)) \inl [(\gamma \vee \alpha) \inl (\gamma \vee \beta)]]$   by axiom $(A{\ref{axioma_menor_cota_superior}})$. \label{280416_15}
			\item 	$\Gamma \cl \gamma \inl (\gamma \vee \beta)$   by axiom $(A{\ref{axioma_supremo_izquierda}})$.  \label{280416_16}
			\item 	$\Gamma \cl (\alpha \inl (\gamma \vee \beta)) \inl [(\gamma \vee \alpha) \inl (\gamma \vee \beta)]$  by ($\N$-MP) applied to \ref{280416_15} and \ref{280416_16}. \label{280416_19}
			\item 	$\Gamma \cl \alpha \inl \beta$ by hypothesis.  \label{280416_17} \label{280416_21}
			\item 	$\Gamma \cl \beta \inl (\gamma \vee \beta)$   by axiom $(A{\ref{axioma_supremo_derecha}})$. \label{280416_18}
			\item 	$\Gamma \cl \alpha \inl (\gamma \vee \beta)$  by axiom $(A{\ref{axioma_transitividadNelson}})$ and ($\N$-MP) applied to \ref{280416_17} and \ref{280416_18}. \label{280416_20}
			\item 	$\Gamma \cl (\gamma \vee \alpha) \inl (\gamma \vee \beta)$  by ($\N$-MP) applied to \ref{280416_19} and \ref{280416_20}.
			\item 	$\Gamma \cl  (\alpha \inl (\beta \vee \gamma)) \inl [(\gamma \inl (\beta \vee \gamma)) \inl [(\alpha \vee \gamma) \inl (\beta \vee \gamma)]]$   by axiom $(A{\ref{axioma_menor_cota_superior}})$. \label{280416_23}
			\item 	$\Gamma \cl \beta \inl (\beta \vee \gamma)$   by axiom $(A{\ref{axioma_supremo_izquierda}})$. \label{280416_22}
			\item 	$\Gamma \cl \alpha \inl (\beta \vee \gamma)$  by axiom $(A{\ref{axioma_transitividadNelson}})$ and ($\N$-MP) applied to \ref{280416_21} and \ref{280416_22}. \label{280416_24}
			\item 	$\Gamma \cl (\gamma \inl (\beta \vee \gamma)) \inl [(\alpha \vee \gamma) \inl (\beta \vee \gamma)]$  by ($\N$-MP) applied to \ref{280416_23} and \ref{280416_24}. \label{280416_25}
			\item 	$\Gamma \cl \gamma \inl (\beta \vee \gamma)$  by axiom $(A{\ref{axioma_supremo_derecha}})$. \label{280416_26}
			\item 	$\Gamma \cl (\alpha \vee \gamma) \inl (\beta \vee \gamma)$  by ($\N$-MP) applied to \ref{280416_25} and \ref{280416_26}.
		\end{enumerate}

		\item[(\ref{080716_42})]
		\begin{enumerate}
			\item 	$\Gamma \cl \alpha \inl (\beta \vee \alpha)$  by axiom $(A{\ref{axioma_supremo_derecha}})$. \label{080716_43}
			\item 	$\Gamma \cl \beta \inl (\beta \vee \alpha)$   by axiom $(A{\ref{axioma_supremo_izquierda}})$. \label{080716_45}
			\item 	$\Gamma \cl (\alpha \inl (\beta \vee \alpha)) \inl [(\beta \inl (\beta \vee \alpha)) \inl ((\alpha \vee \beta) \inl (\beta \vee \alpha))]$  by axiom $(A{\ref{axioma_menor_cota_superior}})$. \label{080716_44}
			\item 	$\Gamma \cl (\beta \inl (\beta \vee \alpha)) \inl ((\alpha \vee \beta) \inl (\beta \vee \alpha))$   by ($\N$-MP) applied to \ref{080716_43} and \ref{080716_44}. \label{080716_46}
			\item 	$\Gamma \cl (\alpha \vee \beta) \inl (\beta \vee \alpha)$   by ($\N$-MP) applied to \ref{080716_45} and \ref{080716_46}.
		\end{enumerate}		
		
		\item[(\ref{290916_00})]
		\begin{enumerate}
			\item 	$\Gamma \cl  (\alpha\land \beta )\inl \alpha$  by axiom $(A{\ref{axioma_infimo_izquierda}})$. \label{290916_01}
			\item 	$\Gamma \cl (\alpha\land \beta) \inl \beta$   by axiom $(A{\ref{axioma_infimo_derecha}})$. \label{290916_02}
			\item 	$\Gamma \cl ((\alpha\land \beta)  \inl \beta) \inl [((\alpha\land \beta) \inl \alpha) \inl ((\alpha \land \beta) \inl (\beta \land \alpha))]$  by axiom $(A{\ref{axioma_mayor_cota_inferior}})$. \label{290916_03}
			\item 	$\Gamma \cl ((\alpha\land \beta) \inl \alpha) \inl ((\alpha \land \beta) \inl (\beta \land \alpha))$   by ($\N$-MP) applied to \ref{290916_02} and \ref{290916_03}. \label{290916_04}
			\item 	$\Gamma \cl (\alpha \land \beta) \inl (\beta \land \alpha)$   by ($\N$-MP) applied to \ref{290916_01} and \ref{290916_04}.
		\end{enumerate}

		\item[(\ref{111215_42})] 
		\begin{enumerate}[1.]
			\item $\Gamma, \alpha \brig \beta \cl  \gamma \inl (\beta \vee \gamma)$ by axiom $(A{\ref{axioma_supremo_derecha}})$.   \label{111215_22}
			\item $\Gamma, \alpha \brig \beta \cl \beta \brig (\beta \vee \gamma) $ by part (\ref{120416_01}).  \label{270417_06}
			\item $\Gamma, \alpha \brig \beta \cl  \alpha \brig \beta$ \label{270417_05} \label{270417_07}
			\item $\Gamma, \alpha \brig \beta \cl  \alpha \inl \beta$ by part (\ref{brigAFactores}) applied to \ref{270417_05}. \label{170417_08}
			\item $\Gamma, \alpha \brig \beta \cl  \alpha \brig (\beta \vee \gamma)$ by part (\ref{propiedad_Transitividad}) applied to \ref{270417_07} and \ref{270417_06}. \label{111215_24}
			\item $\Gamma, \alpha \brig \beta \cl  (\alpha \vee \gamma) \inl (\beta \vee \gamma)$  by part (\ref{280416_14}) applied to \ref{170417_08}.  \label{111215_29}
			\item $\Gamma, \alpha \brig \beta \cl \sim (\beta \vee \gamma) \inl\sim \gamma$  by axiom $(A{\ref{axioma_supremo_negado_derecha}})$.\label{111215_27}
			\item $\Gamma, \alpha \brig \beta \cl  \sim(\beta \vee \gamma) \inl \sim\alpha$ by part (\ref{brigAFactores}) applied to \ref{111215_24}.  \label{111215_25}
			\item $\Gamma, \alpha \brig \beta \cl  (\sim(\beta \vee \gamma) \inl \sim\alpha) \inl [[\sim(\beta \vee \gamma) \inl \sim\gamma] \inl [\sim(\beta \vee \gamma) \inl \sim(\alpha \vee \gamma)]]$  by axiom $(A{\ref{axioma_menor_cota_superior_negado}})$. \label{111215_26}
			\item $\Gamma, \alpha \brig \beta \cl  [\sim(\beta \vee \gamma) \inl \sim\gamma] \inl [\sim(\beta \vee \gamma) \inl \sim(\alpha \vee \gamma)]$  by ($\N$-MP) applied to \ref{111215_25} and \ref{111215_26}. \label{111215_28}
			\item $\Gamma, \alpha \brig \beta \cl  \sim(\beta \vee \gamma) \inl \sim(\alpha \vee \gamma)$  by ($\N$-MP) applied to \ref{111215_27} and \ref{111215_28}.  \label{111215_30}
			\item $\Gamma, \alpha \brig \beta \cl  (\alpha \vee \gamma) \brig (\beta \vee \gamma)$  by part (\ref{021115_01})  applied to \ref{111215_29}  and \ref{111215_30}.
		\end{enumerate}
		
		\item[(\ref{111215_43})] 
		\begin{enumerate}[1.]
			\item $\Gamma, \alpha \brig \beta \cl \beta \brig (\gamma \vee \beta) $ by part (\ref{120416_01}).  \label{170417_10}
			\item $\Gamma, \alpha \brig \beta \cl  \alpha \brig \beta$. \label{170417_09}
			\item $\Gamma, \alpha \brig \beta \cl  \alpha \inl \beta$ by part (\ref{brigAFactores}). \label{120516_02}
			\item $\Gamma, \alpha \brig \beta \cl  \alpha \brig (\gamma \vee \beta)$ by part (\ref{propiedad_Transitividad}) applied to \ref{170417_09} and \ref{170417_10}. \label{111215_32}
			\item $\Gamma, \alpha \brig \beta \cl  (\gamma \vee \alpha) \inl (\gamma \vee \beta)$  by part (\ref{280416_14}) and \ref{120516_02}.  \label{111215_36}
			\item $\Gamma, \alpha \brig \beta \cl  \sim(\gamma \vee \beta) \inl \sim\gamma$  by axiom $(A{\ref{axioma_supremo_negado_izquierda}})$.\label{111215_37}
			\item $\Gamma, \alpha \brig \beta \cl  \sim(\gamma \vee \beta) \inl \sim\alpha$ by part (\ref{brigAFactores}) applied to \ref{111215_32}.  \label{111215_38}
			\item $\Gamma, \alpha \brig \beta \cl  (\sim(\gamma \vee \beta) \inl \sim\gamma) \inl [[\sim(\gamma \vee \beta) \inl \sim\alpha] \inl [\sim(\gamma \vee \beta) \inl \sim(\gamma \vee \alpha)]]$  by axiom $(A{\ref{axioma_menor_cota_superior_negado}})$. \label{111215_39}
			\item $\Gamma, \alpha \brig \beta \cl  [\sim(\gamma \vee \beta) \inl \sim\alpha] \inl [\sim(\gamma \vee \beta) \inl \sim(\gamma \vee \alpha)]$  by ($\N$-MP) applied to \ref{111215_37}  and \ref{111215_39}. \label{111215_40}
			\item $\Gamma, \alpha \brig \beta \cl  \sim(\gamma \vee \beta) \inl \sim(\gamma \vee \alpha)$  by ($\N$-MP) applied to \ref{111215_38} and \ref{111215_40}.  \label{111215_41}
			\item $\Gamma, \alpha \brig \beta \cl  (\gamma \vee \alpha) \brig (\gamma \vee \beta)$  by part (\ref{021115_01})   applied to \ref{111215_36}  and \ref{111215_41}.
		\end{enumerate}
		\item[(\ref{condicionIL3b})] 
		\begin{enumerate}[1.]
			\item $\Gamma, \alpha \brig \beta,  \gamma \brig t \cl  (\alpha \vee \gamma) \brig (\beta \vee \gamma)$ by part (\ref{111215_42}).
			\item $\Gamma, \alpha \brig \beta,  \gamma \brig t  \cl  (\beta \vee \gamma) \brig (\beta \vee t)$ by part (\ref{111215_43}).
			\item $\Gamma, \alpha \brig \beta,  \gamma \brig t  \cl  (\alpha \vee \gamma) \brig (\beta \vee t)$ by part (\ref{propiedad_Transitividad}).
		\end{enumerate}
		
		\item[(\ref{condicionIL3c})] 
		\begin{enumerate}[1.]
			\item $\Gamma, \beta \brig \alpha \cl  \alpha \brig \sim \sim \alpha$ by axiom $(A{\ref{axioma_doble_neg_a_derecha}})$. \label{170417_12}
			\item $\Gamma, \beta \brig \alpha \cl  \beta \brig \alpha$ \label{111215_45} \label{170417_11}
			\item $\Gamma, \beta \brig \alpha \cl \beta \brig \sim \sim \alpha$ by part (\ref{propiedad_Transitividad}) applied to \ref{170417_11} and \ref{170417_12}. \label{170417_14}
			\item $\Gamma, \beta \brig \alpha \cl \sim \sim \beta \brig \beta$ by axiom $(A{\ref{axioma_doble_neg_a_izquierda}})$. \label{170417_13}
			\item $\Gamma, \beta \brig \alpha \cl  \sim \sim \beta \brig \sim \sim \alpha$  by part (\ref{propiedad_Transitividad})  applied to \ref{170417_13} and \ref{170417_14}.
			\item $\Gamma, \beta \brig \alpha \cl  \sim \sim \beta \inl \sim \sim \alpha$ by part (\ref{brigAFactores}).\label{100516_05}
			\item $\Gamma, \beta \brig \alpha \cl \sim \alpha \inl \sim \beta $  by part (\ref{brigAFactores}) applied to \ref{111215_45}. \label{100516_06}
			\item $\Gamma, \beta \brig \alpha \cl  (\sim \alpha \inl \sim \beta) \wedge (\sim \sim \beta \inl \sim \sim \alpha)$ by part (\ref{021115_01}) applied to \ref{100516_05} and \ref{100516_06}. \label{170417_15}
			\item $\Gamma, \beta \brig \alpha \cl \sim \alpha \brig \sim \beta$ by the definition of $\brig$ applied to \ref{170417_15}.
		\end{enumerate}
		
		\item[(\ref{151116_11})]
		\begin{enumerate}
			\item 	$\Gamma \cl (\sim (\alpha \il \beta)) \inl (\alpha \wedge (\sim \beta))$  by axiom $(A{\ref{axioma_paraSN1}})$. \label{151116_12}
			\item 	$\Gamma \cl \sim \beta \inl (\sim (\alpha \wedge \beta))$ by part  (\ref{110716_03}). \label{291116_00}
			\item 	$\Gamma \cl (\alpha \wedge \sim \beta) \inl [\alpha \wedge (\sim (\alpha \wedge \beta))]$ by \ref{291116_00} and part (\ref{280416_01}). \label{151116_13}
			\item 	$\Gamma \cl (\sim (\alpha \il \beta)) \inl [\alpha \wedge (\sim (\alpha \wedge \beta))]$   by axiom $(A\ref{axioma_transitividadNelson})$ and ($\N$-MP) applied to \ref{151116_12} and \ref{151116_13}. \label{151116_14}
			\item 	$\Gamma \cl [\alpha \wedge (\sim (\alpha \wedge \beta))] \inl (\sim (\alpha \inl \beta))$   by axiom $(A{\ref{axioma_paraSN2}})$. \label{151116_15}
			\item 	$\Gamma \cl (\sim (\alpha \il \beta)) \inl (\sim (\alpha \inl \beta))$   by axiom $(A\ref{axioma_transitividadNelson})$ and ($\N$-MP) applied to \ref{151116_14} and \ref{151116_15}.
		\end{enumerate}
		
	\end{itemize}
\end{Proof}


\begin{theorem}(Deduction Theorem) \label{teorema_deduccion}
Let $\Gamma \cup \{\alpha, \beta\} \subseteq \formulas{L}$. Then
	$$\Gamma \cl \alpha \inl \beta \mbox{ if and only if } \Gamma, \alpha \cl \beta $$
\end{theorem}

\begin{Proof} For one implication we have: 
		\begin{enumerate}[1.]
			\item $ \Gamma \cl \alpha \inl \beta$ by hypothesis.
			\item $ \Gamma, \alpha \cl \alpha \inl \beta$. \label{100224_01}
			\item $ \Gamma, \alpha \cl \alpha$. \label{100224_02}
			\item $ \Gamma,\alpha \cl \beta$  by ($\N$-MP) applied to \ref{100224_01} and \ref{100224_02}.
		\end{enumerate}
	For the other one, assume that $\Gamma, \alpha \cl \beta$. We prove the result by induction on the lenght of the proof of $\Gamma, \alpha \cl \beta$.
	\begin{itemize}
		\item If $\cl \beta$ or $\beta \in \Gamma$ then $ \Gamma \cl \beta$ . By Lemma \ref{propiedades_Calculo} (\ref{implicacionATeorema}) we have that $\cl \alpha \inl \beta$. Consequently, $\Gamma \cl \alpha \inl \beta$.
		\item If $\beta = \alpha$, using Lemma \ref{propiedades_Calculo} (\ref{propiedad_XimplicaX}), $\Gamma \cl \alpha \inl \beta$.
		\item If $\beta$ comes from applying the inference rule then there exist $\gamma \in \formulas{L}$ such that
		$\Gamma,\alpha \cl \gamma $ and $\Gamma,\alpha \cl \gamma \inl \beta$.
		Then  
				\begin{enumerate}[1.]
					\item $ \Gamma \cl \alpha \inl (\gamma \inl \beta)$ by inductive hypothesis. \label{250216_09}
					\item $ \Gamma \cl \alpha \inl \gamma$  by inductive hypothesis. \label{250216_13}
					\item $ \Gamma \cl [\alpha \inl (\gamma \inl \beta)] \brig [(\alpha \wedge \gamma) \inl \beta]$  by axiom $(A{\ref{axioma_InfimoAImplicacionVuelta}})$. \label{250216_07}
					\item $ \Gamma \cl [[\alpha \inl (\gamma \inl \beta)] \brig [(\alpha \wedge \gamma) \inl \beta]] \inl [[\alpha \inl (\gamma \inl \beta)] \inl [(\alpha \wedge \gamma) \inl \beta]]$ by axiom $(A{\ref{axioma_infimo_izquierda}})$. \label{250216_08}
					\item $ \Gamma \cl [\alpha \inl (\gamma \inl \beta)] \inl [(\alpha \wedge \gamma) \inl \beta]$  by ($\N$-MP) applied to \ref{250216_07} and \ref{250216_08}. \label{250216_10}
					\item $ \Gamma \cl (\alpha \wedge \gamma) \inl \beta$  by ($\N$-MP) applied to \ref{250216_09} and \ref{250216_10}. \label{250216_17}
					\item $ \Gamma \cl \alpha \inl \alpha$ by Lemma \ref{propiedades_Calculo} (\ref{propiedad_XimplicaX}). \label{250216_11}
					\item $ \Gamma \cl (\alpha \inl \alpha) \inl [(\alpha \inl \gamma) \inl (\alpha \inl (\alpha \wedge \gamma))]$  by axiom $(A{\ref{axioma_mayor_cota_inferior}})$. \label{250216_12}
					\item $ \Gamma \cl (\alpha \inl \gamma) \inl (\alpha \inl (\alpha \wedge \gamma))$  by ($\N$-MP) applied to \ref{250216_11} and \ref{250216_12}. \label{250216_14}
					\item $ \Gamma \cl \alpha \inl (\alpha \wedge \gamma)$  by ($\N$-MP) applied to \ref{250216_13} and \ref{250216_14}. \label{250216_15}
					\item $ \Gamma \cl (\alpha \inl (\alpha \wedge \gamma)) \inl [[(\alpha \wedge \gamma) \inl \beta] \inl [\alpha \inl \beta]]$  by axiom $(A{\ref{axioma_transitividadNelson}})$. \label{250216_16}
					\item $ \Gamma \cl [(\alpha \wedge \gamma) \inl \beta] \inl [\alpha \inl \beta]$  by ($\N$-MP) applied to \ref{250216_15} and \ref{250216_16}. \label{250216_18}
					\item $ \Gamma \cl \alpha \inl \beta$  by ($\N$-MP) applied to \ref{250216_17} and \ref{250216_18}.
				\end{enumerate}
	\end{itemize}
\end{Proof}

We use    $\alpha \leftrightarrow_N \beta$ as an abbreviation for the formula   $(\alpha \inl \beta) \wedge (\beta \inl \alpha)$.

\begin{lema} \label{propiedades_Calculo_conDeduccion_para_reticulado}
	Let $\Gamma \cup \{\alpha, \beta\} \subseteq \formulas{L}$. In $\SNlog$ the following properties hold:
	\begin{enumerate}[{\rm (a)}]
		\item $\Gamma \cl (\sim \alpha \wedge \sim \beta) \leftrightarrow_N \sim (\alpha \vee \beta)$, \label{090416_47}	

		\item $\Gamma \cl \debaj{(\alpha \wedge (\alpha \inl \beta))}{(\alpha \wedge (\sim \alpha \vee \beta))}$, \label{120716_05}	
		\item If $\Gamma \cl \alpha \leftrightarrow_N \beta$ then  $\Gamma \cl \beta \leftrightarrow_N \alpha$, \label{290916_05}
		\item $\Gamma \cl \alpha \brig (\alpha \wedge (\alpha \vee \beta))$,  \label{260216_28} 
		\item $ \Gamma \cl [\alpha \wedge [(\gamma \wedge \alpha) \vee (\beta \wedge \alpha)]] \brig [\alpha \wedge (\beta \vee \gamma)]$, \label{080716_13}
		\item $\Gamma \cl (\alpha \wedge (\beta \vee \gamma)) \brig [\alpha \wedge ((\gamma \wedge \alpha) \vee (\beta \wedge \alpha))]$, \label{260216_01} 
		\item $\Gamma \cl (\sim\beta \inl \sim\alpha) \inl ((\sim\gamma \inl \sim\alpha) \inl (\sim(\beta \wedge \gamma) \inl \sim\alpha))$, \label{110716_01} 
		\item If $ \Gamma \cl \alpha \brig \beta$ and $ \Gamma \cl \alpha \brig \gamma$ then $ \Gamma \cl \alpha \brig \beta \wedge \gamma$ and $\Gamma \cl \alpha \brig \gamma \wedge \beta$, \label{10022016_10} 
		\item 		 \label{100216_12}  
				 If $ \Gamma \cl \alpha \brig \beta$  then $ \Gamma \cl \alpha \brig (\beta \lor \gamma)$ and $\Gamma \cl \alpha \brig (\gamma\lor \beta)$,  
		\item \label{260216_39} 
		 If $\Gamma \cl \alpha \brig \beta$ then  $\Gamma \cl \alpha \land \gamma \brig \beta$ and $\Gamma \cl \gamma \wedge \alpha \brig \beta$, 
		
	\item $\Gamma, \alpha \brig \beta \cl (\alpha \wedge \gamma) \brig (\beta \wedge \gamma)$,   \label{111215_17}
	
	\item $\Gamma, \alpha \brig \beta \cl (\gamma \wedge \alpha) \brig (\gamma \wedge \beta)$,  \label{111215_18}
	
\item $\Gamma, \alpha \brig \beta, \gamma \brig t \cl (\alpha \wedge \gamma) \brig (\beta \wedge t)$,  \label{condicionIL3a}

\item $\Gamma \cl [\sim (\alpha \inl \beta)] \inl (\alpha \wedge \sim \beta)$, \label{090416_19} 
\item $\Gamma \cl (\alpha \wedge \sim \alpha) \inl \beta$, \label{120716_04}  

\item $\Gamma \cl \debaj{(\alpha \wedge \sim \alpha)}{(\beta \vee \sim \beta)}$, \label{120716_01} 

\item $\Gamma \cl (\alpha \inl \beta) \inl [(\beta \inl \alpha) \inl [(\sim (\beta \il \gamma)) \inl (\sim(\alpha \il \gamma))]]$, \label{120716_02} 

\item $\Gamma \cl (\sim \alpha \inl \sim \beta) \inl [(\sim \beta \inl \sim \alpha) \inl [(\sim (\gamma \il \alpha)) \inl (\sim(\gamma \il \beta))]]$, \label{120716_03}   

\item $\Gamma \cl \alpha \brig \beta, \beta \brig \alpha, \gamma \brig t, t \brig \gamma \cl (\alpha \il \gamma) \brig (\beta \il t)$, \label{condicionIL3d} 

\item $\Gamma, \alpha \cl \beta \brig \alpha$, \label{condicionIL5}  

\item $\Gamma \cl (\alpha \il  \beta) \inl (\alpha \inl \beta)$. \label{151116_01}

	\end{enumerate}
\end{lema}

\begin{Proof} 
	\begin{itemize}
		
		\item[(\ref{090416_47})]
		\begin{enumerate}
			\item 	$\Gamma \cl (\sim \sim \alpha \vee \sim \sim \beta) \brig \sim (\sim \alpha \wedge \sim \beta)$  by axiom $(A{\ref{axioma_distribuye_neg_infimo2}})$. \label{090416_54}
			\item 	$\Gamma \cl \alpha \brig (\sim \sim \alpha)$  by axiom $(A{\ref{axioma_doble_neg_a_derecha}})$. \label{090416_48}
			\item 	$\Gamma, \alpha \brig (\sim \sim \alpha) \cl (\alpha \vee \beta) \brig (\sim \sim \alpha \vee \beta)$   by Lemma \ref{propiedades_Calculo} (\ref{111215_42}).
			\item 	$\Gamma \cl (\alpha \brig (\sim \sim \alpha)) \inl [(\alpha \vee \beta) \brig (\sim \sim \alpha \vee \beta)]$  by Theorem \ref{teorema_deduccion}. \label{090416_49}
			\item 	$\Gamma \cl (\alpha \vee \beta) \brig (\sim \sim \alpha \vee \beta)$  by ($\N$-MP) applied to \ref{090416_48} and \ref{090416_49}. \label{090416_52}
			\item 	$\Gamma \cl \beta \brig (\sim \sim \beta)$  by axiom $(A{\ref{axioma_doble_neg_a_derecha}})$. \label{090416_50}
			\item 	$\Gamma, \beta \brig (\sim \sim \beta) \cl (\sim \sim \alpha \vee \beta) \brig (\sim \sim \alpha \vee \sim \sim \beta)$   by Lemma \ref{propiedades_Calculo} (\ref{111215_43}).
			\item 	$\Gamma \cl (\beta \brig (\sim \sim \beta)) \inl [(\sim \sim \alpha \vee \beta) \brig (\sim \sim \alpha \vee \sim \sim \beta)]$   by Theorem \ref{teorema_deduccion}.  \label{090416_51}
			\item 	$\Gamma \cl (\sim \sim \alpha \vee \beta) \brig (\sim \sim \alpha \vee \sim \sim \beta)$  by ($\N$-MP) applied to \ref{090416_50} and \ref{090416_51}. \label{090416_53}
			\item 	$\Gamma \cl (\alpha \vee \beta) \brig (\sim \sim \alpha \vee \sim \sim \beta)$   by Lemma \ref{propiedades_Calculo} (\ref{propiedad_Transitividad}) applied to \ref{090416_52} and \ref{090416_53}. \label{090416_55}
			\item 	$\Gamma \cl (\alpha \vee \beta) \brig \sim (\sim \alpha \wedge \sim \beta)$  by Lemma \ref{propiedades_Calculo} (\ref{propiedad_Transitividad}) applied to \ref{090416_54} and \ref{090416_55}.
			\item 	$\Gamma \cl [\sim \sim (\sim \alpha \wedge \sim \beta)] \inl (\sim (\alpha \vee \beta))$ by Lemma \ref{propiedades_Calculo} (\ref{brigAFactores}). \label{090416_56}
			\item 	$\Gamma \cl (\sim \alpha \wedge \sim \beta) \inl (\sim \sim (\sim \alpha \wedge \sim \beta))$ by axiom $(A{\ref{axioma_doble_neg_a_derecha}})$  and Lemma \ref{propiedades_Calculo} (\ref{brigAFactores}). \label{090416_57}
			\item 	$\Gamma \cl (\sim \alpha \wedge \sim \beta) \inl (\sim (\alpha \vee \beta))$  by axiom $(A{\ref{axioma_transitividadNelson}})$ and ($\N$-MP) applied to \ref{090416_56} and \ref{090416_57}. \label{090416_60}
			\item 	$\Gamma \cl (\sim (\alpha \vee \beta)) \inl (\sim \alpha)$  by axiom $(A{\ref{axioma_supremo_negado_izquierda}})$. \label{090416_58}
			\item 	$\Gamma \cl (\sim (\alpha \vee \beta)) \inl (\sim \beta)$  by axiom $(A{\ref{axioma_supremo_negado_derecha}})$. \label{090416_59}
			\item 	$\Gamma \cl (\sim (\alpha \vee \beta)) \inl (\sim \alpha \wedge \sim \beta)$  by axiom $(A{\ref{axioma_mayor_cota_inferior}})$ and ($\N$-MP) applied to \ref{090416_58} and \ref{090416_59}. \label{090416_61}
			\item 	$\Gamma \cl (\sim \alpha \wedge \sim \beta) \leftrightarrow_N (\sim (\alpha \vee \beta))$  by Lemma \ref{propiedades_Calculo} (\ref{021115_01}) applied to \ref{090416_60} and \ref{090416_61}.
		\end{enumerate}

		\item[(\ref{120716_05})]
		\begin{enumerate}
			\item 	$\Gamma \cl (\alpha \wedge \sim \beta) \inl \sim \beta$  by axiom $(A{\ref{axioma_infimo_derecha}})$.
			\item 	$\Gamma \cl \sim \beta \inl \sim(\alpha \wedge \beta)$   by  \ref{propiedades_Calculo} (\ref{110716_03}).
			\item 	$\Gamma \cl (\alpha \wedge \sim \beta) \inl (\alpha \wedge \sim(\alpha \wedge \beta))$ by Lemma \ref{propiedades_Calculo} (\ref{280416_01}). \label{280416_28}
			\item 	$\Gamma \cl (\alpha \wedge \sim(\alpha \wedge \beta)) \inl (\sim (\alpha \il (\alpha \wedge \beta)))$  by axiom $(A{\ref{axioma_paraSN2}})$. \label{170417_16}
			\item 	$\Gamma \cl (\alpha \wedge \sim(\alpha \wedge \beta)) \inl (\sim (\alpha \inl \beta))$  by definition of $\inl$ applied to \ref{170417_16}. \label{280416_29}  
			
			\item 	$\Gamma \cl (\alpha \wedge \sim \beta) \inl (\sim (\alpha \inl \beta))$  by axiom $(A{\ref{axioma_transitividadNelson}})$ and ($\N$-MP) applied to \ref{280416_28} and \ref{280416_29}.
			\item 	$\Gamma \cl (\sim \alpha \vee (\alpha \wedge \sim \beta)) \inl (\sim \alpha \vee (\sim (\alpha \inl \beta)))$  by Lemma \ref{propiedades_Calculo} (\ref{280416_14}).  \label{280416_38}
			\item 	$\Gamma \cl \sim (\alpha \wedge (\sim \alpha \vee \beta)) \inl (\sim \alpha \vee \sim (\sim \alpha \vee \beta))$ by axiom $(A{\ref{axioma_distribuye_neg_infimo1}})$ and by Lemma \ref{propiedades_Calculo} (\ref{brigAFactores}). \label{280416_33}
			\item 	$\Gamma \cl (\sim \sim \alpha \wedge \sim \beta)  \leftrightarrow_N \sim (\sim \alpha \vee \beta))$  by Lemma \ref{propiedades_Calculo_conDeduccion_para_reticulado} (\ref{090416_47}). \label{280416_31}
			\item 	$\Gamma \cl [(\sim \sim \alpha \wedge \sim \beta)  \leftrightarrow_N \sim (\sim \alpha \vee \beta))] \inl [(\sim (\sim \alpha \vee \beta)) \inl (\sim \sim \alpha \wedge \sim \beta)]$   by axiom $(A{\ref{axioma_infimo_derecha}})$. \label{280416_32}
			\item 	$\Gamma \cl (\sim (\sim \alpha \vee \beta)) \inl (\sim \sim \alpha \wedge \sim \beta)$  by ($\N$-MP) applied to \ref{280416_31} and \ref{280416_32}. \label{170417_17}
			\item 	$\Gamma \cl (\sim \alpha \vee (\sim (\sim \alpha \vee \beta))) \inl (\sim \alpha \vee (\sim \sim \alpha \wedge \sim \beta))$    by Lemma \ref{propiedades_Calculo} (\ref{280416_14}) applied to \ref{170417_17}. \label{280416_34}
			\item 	$\Gamma \cl \sim (\alpha \wedge (\sim \alpha \vee \beta)) \inl (\sim \alpha \vee (\sim \sim \alpha \wedge \sim \beta))$  by axiom $(A{\ref{axioma_transitividadNelson}})$ and ($\N$-MP) applied to \ref{280416_33} and \ref{280416_34}. \label{280416_35}
			\item 	$\Gamma \cl \sim \sim \alpha \inl \alpha$  by axiom $(A{\ref{axioma_doble_neg_a_izquierda}})$ and by Lemma \ref{propiedades_Calculo} (\ref{brigAFactores}). \label{170417_18}
			\item 	$\Gamma \cl (\sim \sim \alpha \wedge \sim \beta) \inl (\alpha  \wedge \sim \beta)$  by Lemma \ref{propiedades_Calculo} (\ref{280416_01}) applied to \ref{170417_18}. \label{170417_19}
			\item 	$\Gamma \cl (\sim \alpha \vee (\sim \sim \alpha \wedge \sim \beta)) \inl (\sim \alpha \vee (\alpha  \wedge \sim \beta))$  by Lemma \ref{propiedades_Calculo} (\ref{280416_14}) applied to \ref{170417_19}. \label{280416_36}
			\item 	$\Gamma \cl (\sim \alpha \vee (\sim (\alpha \inl \beta))) \inl (\sim (\alpha \wedge (\alpha \inl \beta)))$  by axiom $(A{\ref{axioma_distribuye_neg_infimo2}})$ and by Lemma \ref{propiedades_Calculo} (\ref{brigAFactores}). \label{280416_40}
			\item 	$\Gamma \cl \sim (\alpha \wedge (\sim \alpha \vee \beta)) \inl (\sim \alpha \vee (\alpha  \wedge \sim \beta))$  by axiom $(A{\ref{axioma_transitividadNelson}})$ and ($\N$-MP) applied to \ref{280416_35} and \ref{280416_36}. \label{280416_37}
			\item 	$\Gamma \cl \sim (\alpha \wedge (\sim \alpha \vee \beta)) \inl  (\sim \alpha \vee (\sim (\alpha \inl \beta)))$  by axiom $(A{\ref{axioma_transitividadNelson}})$ and ($\N$-MP) applied to \ref{280416_37} and \ref{280416_38}. \label{280416_39}
			\item 	$\Gamma \cl \sim (\alpha \wedge (\sim \alpha \vee \beta)) \inl  (\sim (\alpha \wedge (\alpha \inl \beta)))$  by axiom $(A{\ref{axioma_transitividadNelson}})$ and ($\N$-MP) applied to \ref{280416_39} and \ref{280416_40}. \label{280416_47}
			\item 	$\Gamma, \alpha, \alpha \inl \beta \cl \alpha$. \label{280416_41} \label{170417_20}
			\item 	$\Gamma, \alpha, \alpha \inl \beta \cl \alpha \inl \beta$. \label{280416_42}
			\item 	$\Gamma, \alpha, \alpha \inl \beta \cl \beta$  by ($\N$-MP) applied to \ref{280416_41} and \ref{280416_42}. \label{280416_43}
			\item 	$\Gamma, \alpha, \alpha \inl \beta  \cl \beta \inl (\sim \alpha \vee \beta)$  by axiom $(A{\ref{axioma_supremo_derecha}})$. \label{280416_44}
			\item 	$\Gamma, \alpha, \alpha \inl \beta \cl \sim \alpha \vee \beta$  by ($\N$-MP) applied to \ref{280416_43} and \ref{280416_44}. \label{170417_21}
			\item 	$\Gamma, \alpha, \alpha \inl \beta \cl \alpha \wedge (\sim \alpha \vee \beta)$   by Lemma \ref{propiedades_Calculo} (\ref{021115_01}) applied to \ref{170417_20} and \ref{170417_21}. 
			\item 	$\Gamma, \alpha \cl (\alpha \inl \beta) \inl (\alpha \wedge (\sim \alpha \vee \beta))$   by Theorem \ref{teorema_deduccion}.
			\item 	$\Gamma \cl \alpha \inl [(\alpha \inl \beta) \inl (\alpha \wedge (\sim \alpha \vee \beta))]$   by Theorem \ref{teorema_deduccion}. \label{280416_45}
			\item 	$\Gamma \cl [\alpha \inl [(\alpha \inl \beta) \inl (\alpha \wedge (\sim \alpha \vee \beta))]] \brig [(\alpha \wedge (\alpha \inl \beta)) \inl (\alpha \wedge (\sim \alpha \vee \beta))]$  by axiom $(A{\ref{axioma_InfimoAImplicacionVuelta}})$.
			\item 	$\Gamma \cl [\alpha \inl [(\alpha \inl \beta) \inl (\alpha \wedge (\sim \alpha \vee \beta))]] \inl [(\alpha \wedge (\alpha \inl \beta)) \inl (\alpha \wedge (\sim \alpha \vee \beta))]$ by Lemma \ref{propiedades_Calculo} (\ref{brigAFactores}). \label{280416_46}
			\item 	$\Gamma \cl (\alpha \wedge (\alpha \inl \beta)) \inl (\alpha \wedge (\sim \alpha \vee \beta))$  by ($\N$-MP) applied to \ref{280416_45} and \ref{280416_46}. \label{280416_48}
			\item 	$\Gamma \cl (\alpha \wedge (\alpha \inl \beta)) \brig (\alpha \wedge (\sim \alpha \vee \beta))$    by Lemma \ref{propiedades_Calculo} (\ref{021115_01}) applied to \ref{280416_47} and \ref{280416_48}.
		\end{enumerate}
		\item[(\ref{290916_05})] 
		\begin{enumerate}[1.]
			\item $ \Gamma \cl (\alpha \inl \beta)\land(\beta\inl \alpha)$ by  hypothesis.
			\item $ \Gamma \cl [(\alpha \inl \beta)\land(\beta\inl \alpha)]\inl [(\beta \inl \alpha)\land(\alpha\inl \beta)]$    by  Lemma \ref{propiedades_Calculo} (\ref{290916_00}).
			\item $\Gamma \cl (\beta \inl \alpha)\land (\alpha\inl \beta)$  by ($\N$-MP). 
			\item $ \Gamma \cl \beta \leftrightarrow_N \alpha$ by the definition of $\leftrightarrow_N$.
		\end{enumerate}
		
		\item[ (\ref{260216_28})] 
		\begin{enumerate}[1.]
			\item $ \Gamma \cl \alpha \inl \alpha$ by  Lemma \ref{propiedades_Calculo} (\ref{propiedad_XimplicaX}). \label{260216_29}
			\item $ \Gamma \cl \alpha \inl (\alpha \vee \beta)$   by axiom  $(A{\ref{axioma_supremo_izquierda}})$. \label{260216_30}
			\item $ \Gamma \cl \alpha \inl (\alpha \wedge (\alpha \vee \beta))$  by axiom $(A{\ref{axioma_mayor_cota_inferior}})$ and ($\N$-MP) applied to \ref{260216_29} and \ref{260216_30}. \label{080716_11}
			\item $ \Gamma \cl [\sim (\alpha \wedge (\alpha \vee \beta))] \inl [\sim \alpha \vee \sim (\alpha \vee \beta)]$   by axiom  $(A{\ref{axioma_distribuye_neg_infimo1}})$ and Lemma \ref{propiedades_Calculo} (\ref{brigAFactores}). \label{080716_07}
			\item $ \Gamma \cl \sim (\alpha \vee \beta) \leftrightarrow_N 	(\sim \alpha \wedge \sim \beta)$ by Lemma \ref{propiedades_Calculo_conDeduccion_para_reticulado} (\ref{090416_47}) and (\ref{290916_05}). \label{080716_01}
			\item $ \Gamma \cl [\sim (\alpha \vee \beta) \leftrightarrow_N (\sim \alpha \wedge \sim \beta)] \inl [\sim (\alpha \vee \beta) \inl  (\sim \alpha \wedge \sim \beta)] $  by axiom $(A{\ref{axioma_infimo_izquierda}})$. \label{080716_02}
			\item $ \Gamma \cl \sim (\alpha \vee \beta) \inl  (\sim \alpha \wedge \sim \beta)$   by ($\N$-MP) applied to \ref{080716_01} and \ref{080716_02}. 
			\item $ \Gamma \cl (\sim \alpha \vee \sim (\alpha \vee \beta)) \inl  (\sim \alpha  \vee (\sim \alpha \wedge \sim \beta))$   by Lemma \ref{propiedades_Calculo} (\ref{280416_14}). \label{080716_08} 
			\item $ \Gamma \cl (\sim \alpha \wedge \sim \beta) \inl \sim \alpha$     by axiom $(A{\ref{axioma_infimo_izquierda}})$. \label{080716_05}
			\item $ \Gamma \cl \sim \alpha \inl \sim \alpha$    by Lemma \ref{propiedades_Calculo} (\ref{propiedad_XimplicaX}). \label{080716_03}
			\item $ \Gamma \cl [\sim \alpha \inl \sim \alpha] \inl [[(\sim \alpha \wedge \sim \beta) \inl \sim \alpha] \inl [(\sim \alpha  \vee (\sim \alpha \wedge \sim \beta)) \inl \sim \alpha]]$   by axiom $(A{\ref{axioma_menor_cota_superior}})$. \label{080716_04}
			\item $ \Gamma \cl [(\sim \alpha \wedge \sim \beta) \inl \sim \alpha] \inl [(\sim \alpha  \vee (\sim \alpha \wedge \sim \beta)) \inl \sim \alpha]$  by ($\N$-MP) applied to \ref{080716_03} and \ref{080716_04}. \label{080716_06}
			\item $ \Gamma \cl (\sim \alpha  \vee (\sim \alpha \wedge \sim \beta)) \inl \sim \alpha$   by ($\N$-MP) applied to \ref{080716_05} and \ref{080716_06}. \label{080716_09}
			\item $ \Gamma \cl [\sim (\alpha \wedge (\alpha \vee \beta))] \inl (\sim \alpha  \vee (\sim \alpha \wedge \sim \beta))$   by axiom $(A{\ref{axioma_transitividadNelson}})$ and ($\N$-MP) applied to \ref{080716_07} and \ref{080716_08}. \label{080716_10}
			\item $ \Gamma\cl [\sim (\alpha \wedge (\alpha \vee \beta))]  \inl (\sim \alpha)$  by axiom $(A{\ref{axioma_transitividadNelson}})$ and ($\N$-MP) applied to \ref{080716_09} and \ref{080716_10}. \label{080716_12}
			\item $ \Gamma \cl \alpha \brig (\alpha \wedge (\alpha \vee \beta))$  by Lemma \ref{propiedades_Calculo} (\ref{021115_01}) applied to \ref{080716_11} and \ref{080716_12}.
		\end{enumerate}

		\item[(\ref{080716_13})] 
		\begin{enumerate}[1.]
			\item $ \Gamma \cl (\gamma \wedge \alpha) \inl \alpha$   by axiom $(A{\ref{axioma_infimo_derecha}})$. \label{260216_16}
			\item $ \Gamma \cl (\beta \wedge \alpha) \inl \alpha$  by axiom $(A{\ref{axioma_infimo_derecha}})$. \label{260216_17}
			\item $ \Gamma \cl [(\gamma \wedge \alpha) \vee (\beta \wedge \alpha)] \inl \alpha$  by axiom $(A{\ref{axioma_menor_cota_superior}})$ and ($\N$-MP) applied to \ref{260216_16} and \ref{260216_17}. \label{260216_24}
			\item $ \Gamma \cl (\gamma \wedge \alpha) \inl \gamma$  by axiom $(A{\ref{axioma_infimo_izquierda}})$. \label{260216_18}
			\item $ \Gamma \cl \gamma \inl (\beta \vee \gamma)$ by axiom $(A{\ref{axioma_supremo_derecha}})$. \label{260216_19}
			\item $ \Gamma \cl (\gamma \wedge \alpha) \inl (\beta \vee \gamma)$   by axiom $(A{\ref{axioma_transitividadNelson}})$ and ($\N$-MP) applied to \ref{260216_18} and \ref{260216_19}. \label{260216_22}
			\item $ \Gamma \cl (\beta \wedge \alpha) \inl \beta$ by axiom $(A{\ref{axioma_infimo_izquierda}})$. \label{260216_20}
			\item $ \Gamma \cl \beta \inl (\beta \vee \gamma)$ by axiom $(A{\ref{axioma_supremo_izquierda}})$. \label{260216_21}
			\item $ \Gamma \cl (\beta \wedge \alpha) \inl (\beta \vee \gamma)$  by axiom $(A{\ref{axioma_transitividadNelson}})$ and ($\N$-MP) applied to \ref{260216_20} and \ref{260216_21}. \label{260216_23}
			\item $ \Gamma \cl [(\gamma \wedge \alpha) \vee (\beta \wedge \alpha)] \inl (\beta \vee \gamma)$  by axiom $(A{\ref{axioma_menor_cota_superior}})$ and ($\N$-MP) applied to \ref{260216_22} and \ref{260216_23}. \label{260216_25}
			\item $ \Gamma \cl [(\gamma \wedge \alpha) \vee (\beta \wedge \alpha)] \inl [\alpha \wedge (\beta \vee \gamma)]$  by axiom $(A{\ref{axioma_mayor_cota_inferior}})$ and ($\N$-MP) applied to \ref{260216_24} and \ref{260216_25}. \label{260216_26}
			\item $ \Gamma \cl [\alpha \wedge [(\gamma \wedge \alpha) \vee (\beta \wedge \alpha)]] \inl [(\gamma \wedge \alpha) \vee (\beta \wedge \alpha)]$   by axiom $(A{\ref{axioma_infimo_derecha}})$. \label{260216_27}
			\item $ \Gamma \cl [\alpha \wedge [(\gamma \wedge \alpha) \vee (\beta \wedge \alpha)]] \inl [\alpha \wedge (\beta \vee \gamma)]$  by axiom  $(A{\ref{axioma_transitividadNelson}})$ and ($\N$-MP) applied to \ref{260216_26} and \ref{260216_27}. \label{080716_36}
			\item $ \Gamma \cl \sim (\alpha \wedge (\beta \vee \gamma)) \inl [\sim \alpha \vee \sim (\beta \vee \gamma)]$  by axiom  $(A{\ref{axioma_distribuye_neg_infimo1}})$   and Lemma \ref{propiedades_Calculo} (\ref{brigAFactores}). \label{080716_28}
			\item $ \Gamma \cl \sim (\beta \vee \gamma) \inl \sim \gamma$  by axiom  $(A{\ref{axioma_supremo_negado_derecha}})$. \label{170417_22}
			\item $ \Gamma \cl (\sim \alpha \vee \sim (\beta \vee \gamma)) \inl (\sim \alpha \vee \sim \gamma)$   by Lemma \ref{propiedades_Calculo} (\ref{280416_14}) applied to \ref{170417_22}. \label{080716_14}
			\item $ \Gamma \cl (\sim \alpha \vee \sim \gamma) \inl (\sim \gamma \vee \sim \alpha)$  by Lemma \ref{propiedades_Calculo} (\ref{080716_42}). \label{080716_15} 
			\item $ \Gamma \cl (\sim \alpha \vee \sim (\beta \vee \gamma)) \inl (\sim \gamma \vee \sim \alpha)$  by axiom $(A{\ref{axioma_transitividadNelson}})$ and ($\N$-MP) applied to \ref{080716_14} and \ref{080716_15}. \label{080716_16}
			\item $ \Gamma \cl \sim (\beta \vee \gamma) \inl \sim \beta$  by axiom  $(A{\ref{axioma_supremo_negado_izquierda}})$.
			
			In a similar manner,
			\item $ \Gamma \cl  (\sim \alpha \vee \sim (\beta \vee \gamma)) \inl (\sim \beta \vee \sim \alpha)$. \label{080716_17}
			\item $ \Gamma \cl  (\sim \alpha \vee \sim (\beta \vee \gamma)) \inl [(\sim \gamma \vee \sim \alpha) \wedge (\sim \beta \vee \sim \alpha)]$  by axiom $(A{\ref{axioma_mayor_cota_inferior}})$ and ($\N$-MP) applied to \ref{080716_16} and \ref{080716_17}. \label{080716_29}
			\item $ \Gamma \cl [(\sim \gamma \vee \sim \alpha) \wedge (\sim \beta \vee \sim \alpha)] \inl (\sim \gamma \vee \sim \alpha)$  by axiom $(A{\ref{axioma_infimo_izquierda}})$. \label{080716_18}
			\item $ \Gamma \cl (\sim \gamma \vee \sim \alpha) \inl [\sim (\gamma \wedge \alpha)]$  by axiom $(A{\ref{axioma_distribuye_neg_infimo2}})$   and Lemma \ref{propiedades_Calculo} (\ref{brigAFactores}).  \label{080716_19}
			\item $ \Gamma \cl  [(\sim \gamma \vee \sim \alpha) \wedge (\sim \beta \vee \sim \alpha)] \inl [\sim (\gamma \wedge \alpha)]$  by  axiom $(A{\ref{axioma_transitividadNelson}})$ and ($\N$-MP) applied to \ref{080716_18} and \ref{080716_19}. \label{080716_22}
			\item $ \Gamma \cl [(\sim \gamma \vee \sim \alpha) \wedge (\sim \beta \vee \sim \alpha)] \inl (\sim \beta \vee \sim \alpha)$  by axiom $(A{\ref{axioma_infimo_derecha}})$. \label{080716_20}
			\item $ \Gamma \cl (\sim \beta \vee \sim \alpha) \inl [\sim (\beta \wedge \alpha)]$  by axiom $(A{\ref{axioma_distribuye_neg_infimo2}})$.  \label{080716_21}
			\item $ \Gamma \cl  [(\sim \gamma \vee \sim \alpha) \wedge (\sim \beta \vee \sim \alpha)] \inl [\sim (\beta \wedge \alpha)]$  by  axiom $(A{\ref{axioma_transitividadNelson}})$ and ($\N$-MP) applied to \ref{080716_20} and \ref{080716_21}. \label{080716_23}
			\item $ \Gamma \cl [(\sim \gamma \vee \sim \alpha) \wedge (\sim \beta \vee \sim \alpha)] \inl [[\sim (\gamma \wedge \alpha)] \wedge [\sim (\beta \wedge \alpha)]]$  by axiom $(A{\ref{axioma_mayor_cota_inferior}})$ and ($\N$-MP) applied to \ref{080716_22} and \ref{080716_23}. \label{080716_30}
			
			\item $ \Gamma \cl [\sim [(\gamma \wedge \alpha) \vee (\beta \wedge \alpha)]]  \leftrightarrow_N [(\sim (\gamma \wedge \alpha)) \wedge (\sim (\beta \wedge \alpha))]$ by Lemma \ref{propiedades_Calculo_conDeduccion_para_reticulado} (\ref{090416_47}). \label{080716_24}
			\item $ \Gamma \cl \{[\sim [(\gamma \wedge \alpha) \vee (\beta \wedge \alpha)]]  \leftrightarrow_N [(\sim (\gamma \wedge \alpha)) \wedge (\sim (\beta \wedge \alpha))]\} \inl \{[(\sim (\gamma \wedge \alpha)) \wedge (\sim (\beta \wedge \alpha))] \inl [\sim [(\gamma \wedge \alpha) \vee (\beta \wedge \alpha)]]\}$  by axiom $(A{\ref{axioma_infimo_derecha}})$. \label{080716_25}
			\item $ \Gamma \cl [(\sim (\gamma \wedge \alpha)) \wedge (\sim (\beta \wedge \alpha))] \inl [\sim [(\gamma \wedge \alpha) \vee (\beta \wedge \alpha)]]$  by ($\N$-MP) applied to \ref{080716_24} and \ref{080716_25}. \label{080716_32}
			\item $ \Gamma \cl [\sim [(\gamma \wedge \alpha) \vee (\beta \wedge \alpha)]] \inl [\sim [\alpha \wedge ((\gamma \wedge \alpha) \vee (\beta \wedge \alpha))]]$  by Lemma \ref{propiedades_Calculo} (\ref{110716_03}).  \label{080716_34}
			\item $ \Gamma \cl \sim (\alpha \wedge (\beta \vee \gamma)) \inl [(\sim \gamma \vee \sim \alpha) \wedge (\sim \beta \vee \sim \alpha)]$  by  axiom $(A{\ref{axioma_transitividadNelson}})$ and ($\N$-MP) applied to \ref{080716_28} and \ref{080716_29}. \label{080716_31}
			\item $ \Gamma \cl \sim (\alpha \wedge (\beta \vee \gamma)) \inl [[\sim (\gamma \wedge \alpha)] \wedge [\sim (\beta \wedge \alpha)]]$  by  axiom $(A{\ref{axioma_transitividadNelson}})$ and ($\N$-MP) applied to \ref{080716_30} and \ref{080716_31}. \label{080716_33}
			\item $ \Gamma \cl \sim (\alpha \wedge (\beta \vee \gamma)) \inl  [\sim [(\gamma \wedge \alpha) \vee (\beta \wedge \alpha)]]$  by  axiom $(A{\ref{axioma_transitividadNelson}})$ and ($\N$-MP) applied to \ref{080716_32} and \ref{080716_33}. \label{080716_35}
			\item $ \Gamma \cl \sim (\alpha \wedge (\beta \vee \gamma)) \inl [\sim [\alpha \wedge ((\gamma \wedge \alpha) \vee (\beta \wedge \alpha))]]$  by  axiom $(A{\ref{axioma_transitividadNelson}})$ and ($\N$-MP) applied to \ref{080716_34} and \ref{080716_35}. \label{080716_37}
			\item $ \Gamma \cl [\alpha \wedge [(\gamma \wedge \alpha) \vee (\beta \wedge \alpha)]] \brig [\alpha \wedge (\beta \vee \gamma)]$  by Lemma \ref{propiedades_Calculo} (\ref{021115_01}) applied to \ref{080716_36} and \ref{080716_37}.
		\end{enumerate}

		\item[(\ref{260216_01})] 
		\begin{enumerate}[1.]
			\item $ \Gamma, \alpha, \beta \cl \beta$ \label{170417_23}
			\item $ \Gamma, \beta \cl \alpha \inl \beta$ by Theorem \ref{teorema_deduccion} applied to \ref{170417_23}. \label{260216_04}
			\item $ \Gamma, \beta \cl \alpha \inl \alpha$ by  Lemma \ref{propiedades_Calculo} (\ref{propiedad_XimplicaX}). \label{260216_02}
			\item $ \Gamma, \beta \cl (\alpha \inl \beta) \inl [(\alpha \inl \alpha) \inl (\alpha \inl (\beta \wedge \alpha))]$ by axiom  $(A{\ref{axioma_mayor_cota_inferior}})$. \label{260216_03}
			\item $ \Gamma, \beta \cl (\alpha \inl \alpha) \inl (\alpha \inl (\beta \wedge \alpha))$  by ($\N$-MP) applied to \ref{260216_04} and \ref{260216_03}. \label{260216_05}
			\item $ \Gamma, \beta \cl \alpha \inl (\beta \wedge \alpha)$  by ($\N$-MP) applied to \ref{260216_02} and \ref{260216_05}.
			\item $ \Gamma, \alpha, \beta \cl \beta \wedge \alpha$ by Theorem \ref{teorema_deduccion}. \label{260216_06} \label{170417_24}
			
			\item $ \Gamma, \alpha, \gamma \cl \gamma \wedge \alpha$ similar to \ref{170417_24}. \label{260216_08}
			\item $ \Gamma, \alpha, \beta \cl (\beta \wedge \alpha) \inl [(\gamma \wedge \alpha) \vee (\beta \wedge \alpha)]$ by axiom  $(A{\ref{axioma_supremo_derecha}})$. \label{260216_07}
			\item $ \Gamma, \alpha, \beta \cl (\gamma \wedge \alpha) \vee (\beta \wedge \alpha)$  by ($\N$-MP) applied to \ref{260216_06} and \ref{260216_07}. \label{170417_25}
			\item $ \Gamma, \alpha \cl \beta \inl [(\gamma \wedge \alpha) \vee (\beta \wedge \alpha)]$ by Theorem \ref{teorema_deduccion} applied to \ref{170417_25}. \label{260216_10}
			\item $ \Gamma, \alpha, \gamma \cl (\gamma \wedge \alpha) \inl [(\gamma \wedge \alpha) \vee (\beta \wedge \alpha)]$  by axiom  $(A{\ref{axioma_supremo_izquierda}})$. \label{260216_09}
			\item $ \Gamma, \alpha, \gamma \cl (\gamma \wedge \alpha) \vee (\beta \wedge \alpha)$  by ($\N$-MP) applied to \ref{260216_08} and \ref{260216_09}. \label{170417_26}
			\item $ \Gamma, \alpha \cl  \gamma \inl [(\gamma \wedge \alpha) \vee (\beta \wedge \alpha)]$  by Theorem \ref{teorema_deduccion} applied to \ref{170417_26}. \label{260216_12}
			\item $ \Gamma, \alpha \cl [\beta \inl [(\gamma \wedge \alpha) \vee (\beta \wedge \alpha)]] \inl [[\gamma \inl [(\gamma \wedge \alpha) \vee (\beta \wedge \alpha)]] \inl [(\beta \vee \gamma) \inl [(\gamma \wedge \alpha) \vee (\beta \wedge \alpha)]]]$ by axiom  $(A{\ref{axioma_menor_cota_superior}})$. \label{260216_11}
			\item $ \Gamma, \alpha \cl [\gamma \inl [(\gamma \wedge \alpha) \vee (\beta \wedge \alpha)]] \inl [(\beta \vee \gamma) \inl [(\gamma \wedge \alpha) \vee (\beta \wedge \alpha)]]$  by ($\N$-MP) applied to \ref{260216_10} and \ref{260216_11}. \label{260216_13}
			\item $ \Gamma, \alpha \cl (\beta \vee \gamma) \inl [(\gamma \wedge \alpha) \vee (\beta \wedge \alpha)]$  by ($\N$-MP) applied to \ref{260216_12} and \ref{260216_13}. \label{170417_27}
			\item $ \Gamma \cl \alpha \inl [(\beta \vee \gamma) \inl [(\gamma \wedge \alpha) \vee (\beta \wedge \alpha)]]$  by Theorem \ref{teorema_deduccion} applied to \ref{170417_27}. \label{170417_28}
			\item $ \Gamma \cl [\alpha \inl [(\beta \vee \gamma) \inl [(\gamma \wedge \alpha) \vee (\beta \wedge \alpha)]]] \brig [(\alpha \wedge (\beta \vee \gamma)) \inl [(\gamma \wedge \alpha) \vee (\beta \wedge \alpha)]]$  by axiom  $(A{\ref{axioma_InfimoAImplicacionVuelta}})$.
			\item $ \Gamma \cl (\alpha \wedge (\beta \vee \gamma)) \inl [(\gamma \wedge \alpha) \vee (\beta \wedge \alpha)]$  by Lemma \ref{propiedades_Calculo} (\ref{MPsobreBrigN}) applied to \ref{170417_28}. \label{260216_14}
			\item $ \Gamma \cl (\alpha \wedge (\beta \vee \gamma)) \inl \alpha$  by axiom  $(A{\ref{axioma_infimo_izquierda}})$. \label{260216_15}
			\item $ \Gamma \cl (\alpha \wedge (\beta \vee \gamma))  \inl [\alpha \wedge  [(\gamma \wedge \alpha) \vee (\beta \wedge \alpha)]]$  by axiom $(A{\ref{axioma_mayor_cota_inferior}})$ and ($\N$-MP) applied to \ref{260216_14} and \ref{260216_15}. \label{080716_38}
			\item $ \Gamma \cl [\sim (\alpha \wedge ((\gamma \wedge \alpha) \vee (\beta \wedge \alpha)))] \inl [\sim (\alpha \wedge (\beta \vee \gamma))]$  by axiom  $(A{\ref{axioma_para_reticulado_neg}})$. \label{080716_39}
			\item $\Gamma \cl (\alpha \wedge (\beta \vee \gamma)) \brig [\alpha \wedge ((\gamma \wedge \alpha) \vee (\beta \wedge \alpha))]$  by Lemma \ref{propiedades_Calculo} (\ref{021115_01}) applied to \ref{080716_38} and \ref{080716_39}.
			
		\end{enumerate}
		
		\item[(\ref{110716_01})] 
		\begin{enumerate}
			\item 	$\Gamma, \sim \beta \inl \sim \alpha, \sim \gamma \inl \sim \alpha, \sim(\beta \wedge \gamma) \cl (\sim \beta \inl \sim \alpha) \inl [(\sim \gamma \inl \sim \alpha) \inl ((\sim \beta \vee \sim \gamma) \inl \sim \alpha)]$  by axiom $(A{\ref{axioma_menor_cota_superior}})$. \label{190416_05}
			\item 	$\Gamma, \sim \beta \inl \sim \alpha, \sim \gamma \inl \sim \alpha, \sim(\beta \wedge \gamma) \cl \sim \beta \inl \sim \alpha$.  \label{190416_06}
			\item 	$\Gamma, \sim \beta \inl \sim \alpha, \sim \gamma \inl \sim \alpha, \sim(\beta \wedge \gamma) \cl (\sim \gamma \inl \sim \alpha) \inl ((\sim \beta \vee \sim \gamma) \inl \sim \alpha)$ by ($\N$-MP) applied to \ref{190416_05} and \ref{190416_06}. \label{190416_07}
			\item 	$\Gamma, \sim \beta \inl \sim \alpha, \sim \gamma \inl \sim \alpha, \sim(\beta \wedge \gamma) \cl \sim \gamma \inl \sim \alpha$. \label{190416_08}
			\item 	$\Gamma, \sim \beta \inl \sim \alpha, \sim \gamma \inl \sim \alpha, \sim(\beta \wedge \gamma) \cl (\sim \beta \vee \sim \gamma) \inl \sim \alpha$  by ($\N$-MP) applied to \ref{190416_07} and \ref{190416_08}. \label{190416_11}
			\item 	$\Gamma, \sim \beta \inl \sim \alpha, \sim \gamma \inl \sim \alpha, \sim(\beta \wedge \gamma) \cl \sim(\beta \wedge \gamma)$. \label{190416_09}
			\item 	$\Gamma, \sim \beta \inl \sim \alpha, \sim \gamma \inl \sim \alpha, \sim(\beta \wedge \gamma) \cl [\sim(\beta \wedge \gamma)] \brig (\sim \beta \vee \sim \gamma)$  by axiom $(A{\ref{axioma_distribuye_neg_infimo1}})$.
			\item 	$\Gamma, \sim \beta \inl \sim \alpha, \sim \gamma \inl \sim \alpha, \sim(\beta \wedge \gamma) \cl [\sim(\beta \wedge \gamma)] \inl (\sim \beta \vee \sim \gamma)$  by Lemma \ref{propiedades_Calculo} (\ref{brigAFactores}). \label{190416_10}
			\item 	$\Gamma, \sim \beta \inl \sim \alpha, \sim \gamma \inl \sim \alpha, \sim(\beta \wedge \gamma)  \cl \sim \beta \vee \sim \gamma$  by ($\N$-MP) applied to \ref{190416_09} and \ref{190416_10}. \label{190416_12}
			\item 	$\Gamma, \sim \beta \inl \sim \alpha, \sim \gamma \inl \sim \alpha, \sim(\beta \wedge \gamma)  \cl  \sim \alpha$  by ($\N$-MP) applied to \ref{190416_11} and \ref{190416_12}. \label{170417_29}
			\item 	$\Gamma, \sim \beta \inl \sim \alpha, \sim \gamma \inl \sim \alpha \cl (\sim(\beta \wedge \gamma)) \inl (\sim \alpha)$  by Theorem \ref{teorema_deduccion} applied to \ref{170417_29}. \label{170417_30}
			\item 	$\Gamma, \sim \beta \inl \sim \alpha \cl (\sim \gamma \inl \sim \alpha) \inl [(\sim(\beta \wedge \gamma)) \inl (\sim \alpha)]$  by Theorem \ref{teorema_deduccion} applied to \ref{170417_30}. \label{170417_31}
			\item 	$\Gamma \cl (\sim \beta \inl \sim \alpha) \inl [(\sim \gamma \inl \sim \alpha) \inl [(\sim(\beta \wedge \gamma)) \inl (\sim \alpha)]]$  by Theorem \ref{teorema_deduccion} applied to \ref{170417_31}.
		\end{enumerate}

		\item[ (\ref{10022016_10})] 
		\begin{enumerate}[1.]
			\item $ \Gamma \cl \alpha \brig \beta$ by hypothesis.
			\item $ \Gamma \cl \alpha \inl \beta$  by Lemma \ref{propiedades_Calculo} (\ref{brigAFactores}). \label{100216_01}
			\item $ \Gamma \cl \sim\beta \inl \sim\alpha$   by Lemma \ref{propiedades_Calculo} (\ref{brigAFactores}). \label{100216_05}
			\item $ \Gamma \cl \alpha \brig \gamma$ by hypothesis.
			\item $ \Gamma \cl \alpha \inl \gamma$  by Lemma \ref{propiedades_Calculo} (\ref{brigAFactores}). \label{100216_03}
			\item $ \Gamma \cl \sim\gamma \inl \sim\alpha$   by Lemma \ref{propiedades_Calculo} (\ref{brigAFactores}). \label{100216_07}
			\item $ \Gamma \cl (\alpha \inl \beta) \inl ((\alpha \inl \gamma) \inl (\alpha \inl (\beta \wedge \gamma)))$  by axiom $(A{\ref{axioma_mayor_cota_inferior}})$. \label{100216_02}
			\item $ \Gamma \cl (\alpha \inl \gamma) \inl (\alpha \inl (\beta \wedge \gamma))$  by ($\N$-MP) applied to \ref{100216_01} and \ref{100216_02}. \label{100216_04}
			\item $ \Gamma \cl \alpha \inl (\beta \wedge \gamma)$ \label{100516_01}  by ($\N$-MP) applied to \ref{100216_03} and \ref{100216_04}.
			\item $ \Gamma \cl (\sim\beta \inl \sim\alpha) \inl ((\sim\gamma \inl \sim\alpha) \inl (\sim(\beta \wedge \gamma) \inl \sim\alpha))$  by (\ref{110716_01}).
			\label{100216_06}
			\item $ \Gamma \cl (\sim\gamma \inl \sim\alpha) \inl (\sim(\beta \wedge \gamma) \inl \sim\alpha)$  by ($\N$-MP) applied to \ref{100216_05} and \ref{100216_06}. \label{100216_08}
			\item $ \Gamma \cl \sim(\beta \wedge \gamma) \inl \sim\alpha $  \label{100516_02} by ($\N$-MP) applied to \ref{100216_07} and \ref{100216_08}.
			\item $ \Gamma \cl \alpha \brig \beta \wedge \gamma$ by Lemma \ref{propiedades_Calculo} (\ref{021115_01}) applied to \ref{100516_01} and \ref{100516_02}.
			
			In a similar manner, 
			\item $\Gamma \cl \alpha \brig \gamma \wedge \beta$.
		\end{enumerate}

		\item[(\ref{100216_12})] 
		
		\begin{enumerate}[1.]
			\item $ \Gamma \cl \alpha \brig \beta$ by hypothesis. \label{300916_00}
			\item $ \Gamma \cl \alpha \inl \beta$  by \ref{propiedades_Calculo} (\ref{brigAFactores}).  \label{300916_01}
			\item $\Gamma \cl \beta\inl (\beta\lor \gamma)$  $(A{\ref{axioma_supremo_izquierda}})$. \label{300916_03}
			\item $\Gamma \cl \alpha\inl (\beta\lor \gamma)$  by axiom $(A{\ref{axioma_transitividadNelson}})$ and ($\N$-MP) applied to \ref{300916_01} and \ref{300916_03}. \label{300916_04}
			\item $ \Gamma \cl \sim \beta \inl \sim \alpha$  by \ref{propiedades_Calculo} (\ref{brigAFactores}) applied to \ref{300916_00}. \label{300916_06} 
			\item $ \Gamma \cl \sim(\beta\lor \gamma) \inl \sim \beta$  by axiom $(A{\ref{axioma_supremo_negado_izquierda}})$. \label{300916_07} 
			\item $ \Gamma \cl \sim(\beta\lor \gamma) \inl \sim \alpha$  by axiom $(A{\ref{axioma_transitividadNelson}})$ and ($\N$-MP) applied to \ref{300916_06} and \ref{300916_07}. \label{300916_08} \label{300916_09}  
			\item $ \Gamma \cl \alpha \brig (\beta\lor \gamma)$  by \ref{propiedades_Calculo} (\ref{021115_01}) applied to \ref{300916_04} and \ref{300916_09}.
			
			In a similar manner, using axioms $(A{\ref{axioma_supremo_derecha}})$ and $(A{\ref{axioma_supremo_negado_derecha}})$,
			
			\item $ \Gamma \cl \alpha \brig (\gamma\lor \beta)$  
		\end{enumerate}

		\item[(\ref{260216_39})]  Follows from \ref{propiedades_Calculo}  (\ref{propiedad_Transitividad}) and \ref{propiedades_Calculo}  (\ref{260216_37}).
		
		\item[(\ref{111215_17})] 
		\begin{enumerate}[1.]
			\item $\Gamma, \alpha \brig \beta \cl \alpha \brig \beta$.  \label{021215_01}
			\item $\Gamma, \alpha \brig \beta \cl \alpha \inl \beta$ by  \ref{propiedades_Calculo} (\ref{brigAFactores}).
			\item $\Gamma, \alpha \brig \beta \cl (\alpha \wedge \gamma) \brig \alpha$ by \ref{propiedades_Calculo} (\ref{260216_37}).
			\label{021215_02}
			\item $\Gamma, \alpha \brig \beta \cl (\alpha \wedge \gamma) \brig \beta$ by  \ref{propiedades_Calculo} (\ref{propiedad_Transitividad}) applied to \ref{021215_01} and \ref{021215_02}. \label{021215_07}
			\item $\Gamma, \alpha \brig \beta \cl (\alpha \wedge \gamma) \inl (\beta \wedge \gamma)$  by Lemma \ref{propiedades_Calculo} (\ref{brigAFactores}) and  (\ref{280416_01}). \label{021215_06}
			\item $\Gamma, \alpha \brig \beta \cl \sim\beta \inl \sim(\alpha \wedge \gamma)$ by  \ref{propiedades_Calculo} (\ref{brigAFactores}) applied to \ref{021215_07}. \label{021215_08}
			\item $\Gamma, \alpha \brig \beta \cl \sim\gamma \inl \sim(\alpha \wedge \gamma)$  by   \ref{propiedades_Calculo} (\ref{110716_03}). \label{021215_10}
			\item $\Gamma, \alpha \brig \beta \cl (\sim\beta \inl \sim(\alpha \wedge \gamma)) \inl ((\sim\gamma \inl \sim(\alpha \wedge \gamma)) \inl (\sim(\beta \wedge \gamma) \inl \sim(\alpha \wedge \gamma)))$  by 
			Lema \ref{propiedades_Calculo_conDeduccion_para_reticulado} (\ref{110716_01}).
			\label{021215_09}
			\item $\Gamma, \alpha \brig \beta \cl (\sim\gamma \inl \sim(\alpha \wedge \gamma)) \inl (\sim(\beta \wedge \gamma) \inl \sim(\alpha \wedge \gamma))$  by ($\N$-MP) applied to \ref{021215_08} and \ref{021215_09}. \label{021215_11}
			\item $\Gamma, \alpha \brig \beta \cl \sim(\beta \wedge \gamma) \inl \sim(\alpha \wedge \gamma)$  by ($\N$-MP) applied to \ref{021215_10} and \ref{021215_11}. \label{021215_13}
			\item $\Gamma, \alpha \brig \beta \cl  (\alpha \wedge \gamma) \brig (\beta \wedge \gamma)$ by  \ref{propiedades_Calculo} (\ref{021115_01})  applied to \ref{021215_06} and \ref{021215_13}.
		\end{enumerate}
		
		\item[ (\ref{111215_18})] 
		\begin{enumerate}[1.]
			\item $\Gamma, \alpha \brig \beta \cl (\gamma \wedge \alpha) \brig \alpha $ by  \ref{propiedades_Calculo} (\ref{260216_37}). 
			\label{111215_03}
			\item $\Gamma, \alpha \brig \beta \cl  \alpha \brig \beta$. \label{111215_04}
			\item $\Gamma, \alpha \brig \beta \cl  \alpha \inl \beta$ \label{120516_01} by  \ref{propiedades_Calculo} (\ref{brigAFactores}).
			\item $\Gamma, \alpha \brig \beta \cl  (\gamma \wedge \alpha) \brig \beta$   by  \ref{propiedades_Calculo} (\ref{propiedad_Transitividad})  applied to \ref{111215_03} and \ref{111215_04}.  \label{111215_09}
			\item $\Gamma, \alpha \brig \beta \cl  (\gamma \wedge \alpha) \inl (\gamma \wedge \beta)$  by \ref{120516_01} and Lemma \ref{propiedades_Calculo} (\ref{280416_01}). \label{111215_15} 
			\item $\Gamma, \alpha \brig \beta \cl  \sim\beta \inl \sim(\gamma \wedge \alpha)$ by Lemma \ref{propiedades_Calculo} (\ref{brigAFactores}) applied to \ref{111215_09}. \label{111215_13}
			\item $\Gamma, \alpha \brig \beta \cl  \sim\gamma \inl \sim(\gamma \wedge \alpha)$ by Lemma  \ref{propiedades_Calculo} (\ref{110716_02}). \label{111215_11}
			\item $\Gamma, \alpha \brig \beta \cl  (\sim\gamma \inl \sim(\gamma \wedge \alpha)) \inl [(\sim\beta \inl \sim(\gamma \wedge \alpha)) \inl (\sim(\gamma \wedge \beta) \inl \sim(\gamma \wedge \alpha))]$  by 
			Lemma \ref{propiedades_Calculo_conDeduccion_para_reticulado} (\ref{110716_01}).
			\label{111215_12}
			\item $\Gamma, \alpha \brig \beta \cl  (\sim\beta \inl \sim(\gamma \wedge \alpha)) \inl (\sim(\gamma \wedge \beta) \inl \sim(\gamma \wedge \alpha))$  by ($\N$-MP) applied to \ref{111215_11} and \ref{111215_12}.  \label{111215_14} 
			\item $\Gamma, \alpha \brig \beta \cl  \sim(\gamma \wedge \beta) \inl \sim(\gamma \wedge \alpha)$ by ($\N$-MP) applied to \ref{111215_13} and \ref{111215_14}. \label{111215_16} 
			\item $\Gamma, \alpha \brig \beta \cl  (\gamma \wedge \alpha) \brig (\gamma \wedge \beta)$ by Lemma \ref{propiedades_Calculo} (\ref{021115_01})   applied to \ref{111215_15}  and \ref{111215_16}.
		\end{enumerate}
		
		\item[(\ref{condicionIL3a})] 
		\begin{enumerate}[1.]
			\item $\Gamma, \alpha \brig \beta,  \gamma \brig t   \cl  (\alpha \wedge \gamma) \brig (\beta \wedge \gamma)$ by part (\ref{111215_17}).
			\item $\Gamma, \alpha \brig \beta,  \gamma \brig t \cl  (\beta \wedge \gamma) \brig (\beta \wedge t)$ by part (\ref{111215_18}).
			\item $\Gamma, \alpha \brig \beta,  \gamma \brig t  \cl  (\alpha \wedge \gamma) \brig (\beta \wedge t)$ by Lemma \ref{propiedades_Calculo} (\ref{propiedad_Transitividad}).
		\end{enumerate}
		
		\item[ (\ref{090416_19})]
		\begin{enumerate}
			
			\item 	$\Gamma \cl [\sim (\alpha \il (\alpha \wedge \beta))] \inl [\alpha \land (\sim (\alpha \wedge \beta))]$ by axiom $(A{\ref{axioma_paraSN1}})$ \label{170417_32}
			\item 	$\Gamma \cl [\sim (\alpha \inl \beta)] \inl [\alpha \land (\sim (\alpha \wedge \beta))]$  by definition of $\inl$ applied to  \ref{170417_32}. \label{090416_22}			
			
			\item 	$\Gamma, (\sim (\alpha \wedge \beta)) \brig (\sim \alpha \vee \sim \beta) \cl [\alpha \wedge (\sim (\alpha \wedge \beta))] \brig [\alpha \wedge (\sim \alpha \vee \sim \beta)]$  by  part (\ref{111215_18}). \label{170417_33}
			\item 	$\Gamma \cl [(\sim (\alpha \wedge \beta)) \brig (\sim \alpha \vee \sim \beta)] \inl [[\alpha \wedge (\sim (\alpha \wedge \beta))] \brig [\alpha \wedge (\sim \alpha \vee \sim \beta)]]$  by Theorem \ref{teorema_deduccion} applied to \ref{170417_33}. \label{090416_20}
			\item 	$\Gamma \cl (\sim (\alpha \wedge \beta)) \brig (\sim \alpha \vee \sim \beta)$   by axiom $(A{\ref{axioma_distribuye_neg_infimo1}})$. \label{090416_21}
			\item 	$\Gamma \cl [\alpha \wedge (\sim (\alpha \wedge \beta))] \brig [\alpha \wedge (\sim \alpha \vee \sim \beta)]$   by ($\N$-MP) applied to \ref{090416_20} and \ref{090416_21}.
			\item 	$\Gamma \cl [\alpha \wedge (\sim (\alpha \wedge \beta))] \inl [\alpha \wedge (\sim \alpha \vee \sim \beta)]$ by Lemma \ref{propiedades_Calculo} (\ref{brigAFactores}). \label{090416_23}
			\item 	$\Gamma \cl [\sim (\alpha \inl \beta)] \inl [\alpha \wedge (\sim \alpha \vee \sim \beta)]$  by axiom \label{090416_24} $(A{\ref{axioma_transitividadNelson}})$ and ($\N$-MP) applied to \ref{090416_22} and \ref{090416_23}.
			\item 	$\Gamma \cl [\alpha \wedge (\sim \alpha \vee \sim \beta)] \inl [\alpha \wedge (\alpha \inl (\sim \beta))]$  by axiom $(A{\ref{axioma_implica_infimo_dosVariables2}})$  and Lemma \ref{propiedades_Calculo} (\ref{brigAFactores}). \label{090416_25}
			\item 	$\Gamma \cl [\sim (\alpha \inl \beta)] \inl [\alpha \wedge (\alpha \inl (\sim \beta))]$  by axiom $(A{\ref{axioma_transitividadNelson}})$ and ($\N$-MP) applied to \ref{090416_24} and \ref{090416_25}. \label{090416_26} \label{090416_28}
			\item 	$\Gamma \cl [\alpha \wedge (\alpha \inl (\sim \beta))] \inl \alpha$  by axiom $(A{\ref{axioma_infimo_izquierda}})$. \label{090416_27}
			\item 	$\Gamma \cl [\sim (\alpha \inl \beta)] \inl \alpha$  by axiom $(A{\ref{axioma_transitividadNelson}})$ and ($\N$-MP) applied to \ref{090416_26} and \ref{090416_27}.
			\item 	$\Gamma, \sim (\alpha \inl \beta) \cl \alpha$   by Theorem \ref{teorema_deduccion}. \label{090416_30}
			\item 	$\Gamma \cl [\alpha \wedge (\alpha \inl (\sim \beta))] \inl (\alpha \inl (\sim \beta)) $ by axiom $(A{\ref{axioma_infimo_derecha}})$.  \label{090416_29}
			\item 	$\Gamma \cl [\sim (\alpha \inl \beta)] \inl (\alpha \inl (\sim \beta)) $  by axiom $(A{\ref{axioma_transitividadNelson}})$ and ($\N$-MP) applied to \ref{090416_28} and \ref{090416_29}. \label{170417_34}
			\item 	$\Gamma, \sim (\alpha \inl \beta) \cl \alpha \inl (\sim \beta)$   by Theorem \ref{teorema_deduccion} applied to \ref{170417_34}. \label{090416_31}
			\item 	$\Gamma, \sim (\alpha \inl \beta) \cl \sim \beta$  by ($\N$-MP) applied to \ref{090416_30} and \ref{090416_31}. \label{090416_33}
			\item 	$\Gamma, \sim (\alpha \inl \beta)  \cl (\alpha \wedge \sim \beta)$ by Lemma \ref{propiedades_Calculo} (\ref{021115_01})   applied to \ref{090416_30} and \ref{090416_33}.
			\item 	$\Gamma \cl [\sim (\alpha \inl \beta)] \inl (\alpha \wedge \sim \beta)$   by Theorem \ref{teorema_deduccion}.
		\end{enumerate}
		
		\item[(\ref{120716_04})] 
		\begin{enumerate} 
			\item 	$\Gamma, \alpha, \sim \alpha \cl \sim \alpha$. \label{270416_01}
			\item 	$\Gamma, \alpha, \sim \alpha \cl \sim \alpha \inl (\sim \alpha \vee \beta)$  by axiom $(A{\ref{axioma_supremo_izquierda}})$.  \label{270416_02}
			\item 	$\Gamma, \alpha, \sim \alpha \cl \sim \alpha \vee \beta$  by ($\N$-MP) applied to \ref{270416_01} and \ref{270416_02}. \label{170417_35}
			\item 	$\Gamma, \alpha, \sim \alpha \cl \alpha$. \label{270416_07} \label{170417_36}
			\item 	$\Gamma, \alpha, \sim \alpha  \cl \alpha \wedge (\sim \alpha \vee \beta)$ by Lemma \ref{propiedades_Calculo} (\ref{021115_01}) applied to \ref{170417_35} and \ref{170417_36}. \label{270416_03}
			\item 	$\Gamma, \alpha, \sim \alpha  \cl [\alpha \wedge (\sim \alpha \vee \beta)] \inl [\alpha \wedge (\alpha \inl \beta)]$   by axiom $(A{\ref{axioma_implica_infimo_dosVariables2}})$ and Lemma \ref{propiedades_Calculo} (\ref{brigAFactores}). \label{270416_04}   
			\item 	$\Gamma, \alpha, \sim \alpha  \cl \alpha \wedge (\alpha \inl \beta)$  by ($\N$-MP) applied to \ref{270416_03} and \ref{270416_04}. \label{270416_05}
			\item 	$\Gamma, \alpha, \sim \alpha  \cl (\alpha \wedge (\alpha \inl \beta)) \inl (\alpha \inl \beta)$  by axiom $(A{\ref{axioma_infimo_derecha}})$. \label{270416_06}
			\item 	$\Gamma, \alpha, \sim \alpha  \cl \alpha \inl \beta$  by ($\N$-MP) applied to \ref{270416_05} and \ref{270416_06}. \label{270416_08}
			\item 	$\Gamma, \alpha, \sim \alpha  \cl \beta$  by ($\N$-MP) applied to \ref{270416_07} and \ref{270416_08}. \label{170417_37}
			\item 	$\Gamma, \alpha  \cl \sim \alpha \inl \beta$   by Theorem \ref{teorema_deduccion} applied to \ref{170417_37}. \label{170417_38}
			\item 	$\Gamma \cl \alpha \inl (\sim \alpha \inl \beta)$   by Theorem \ref{teorema_deduccion} applied to \ref{170417_38}. \label{270416_09}
			\item 	$\Gamma \cl [\alpha \inl (\sim \alpha \inl \beta)] \inl [(\alpha \wedge \sim \alpha) \inl \beta]$   by axiom $(A{\ref{axioma_InfimoAImplicacionVuelta}})$ and Lemma \ref{propiedades_Calculo} (\ref{brigAFactores}). \label{270416_10}
			\item 	$\Gamma \cl (\alpha \wedge \sim \alpha) \inl \beta$  by ($\N$-MP) applied to \ref{270416_09} and \ref{270416_10}.
		\end{enumerate}

		\item[(\ref{120716_01})]	
		\begin{enumerate}
			\item 	$\Gamma \cl (\sim \beta \wedge \sim \sim \beta) \leftrightarrow_N \sim(\beta \vee \sim \beta)$  by  \ref{propiedades_Calculo_conDeduccion_para_reticulado} (\ref{090416_47}). \label{210416_11}
			\item 	$\Gamma \cl [(\sim \beta \wedge \sim \sim \beta) \leftrightarrow_N \sim(\beta \vee \sim \beta)] \inl [\sim(\beta \vee \sim \beta) \inl (\sim \beta \wedge \sim \sim \beta)]$ by axiom $(A{\ref{axioma_infimo_derecha}})$. \label{210416_12}
			\item 	$\Gamma \cl \sim(\beta \vee \sim \beta) \inl (\sim \beta \wedge \sim \sim \beta)$  by ($\N$-MP) applied to \ref{210416_11} and \ref{210416_12}. \label{210416_13}
			\item 	$\Gamma \cl (\sim \beta \wedge \sim \sim \beta) \inl \sim (\alpha \wedge \sim \alpha)$   by \ref{propiedades_Calculo_conDeduccion_para_reticulado} (\ref{120716_04}). \label{210416_14}
			\item 	$\Gamma \cl \sim(\beta \vee \sim \beta) \inl \sim (\alpha \wedge \sim \alpha)$  by axiom $(A{\ref{axioma_transitividadNelson}})$ and ($\N$-MP) applied to \ref{210416_13} and \ref{210416_14}. \label{170417_40}
			\item 	$\Gamma \cl (\alpha \wedge \sim \alpha) \inl (\beta \vee \sim \beta)$  by  \ref{propiedades_Calculo_conDeduccion_para_reticulado} (\ref{120716_04}). \label{170417_39}
			\item 	$\Gamma \cl (\alpha \wedge \sim \alpha) \brig (\beta \vee \sim \beta)$  by Lemma \ref{propiedades_Calculo} (\ref{021115_01}) applied to \ref{170417_39} and \ref{170417_40}.
		\end{enumerate}

		\item[(\ref{120716_02})]
		\begin{enumerate}
			\item 	$\Gamma, \alpha \inl \beta, \beta \inl \alpha, \beta \cl \beta \inl \alpha$. \label{210416_15}
			\item 	$\Gamma, \alpha \inl \beta, \beta \inl \alpha, \beta \cl \beta$. \label{210416_16}
			\item 	$\Gamma, \alpha \inl \beta, \beta \inl \alpha, \beta \cl \alpha$  by ($\N$-MP) applied to \ref{210416_15} and \ref{210416_16}. \label{170417_41}
			\item 	$\Gamma, \alpha \inl \beta, \beta \inl \alpha, \beta \cl \sim \gamma \inl \alpha$  by Lemma \ref{propiedades_Calculo} (\ref{implicacionATeorema}) applied to \ref{170417_41}. \label{170417_42}
			\item 	$\Gamma, \alpha \inl \beta, \beta \inl \alpha \cl \beta \inl (\sim \gamma \inl \alpha)$  by Theorem \ref{teorema_deduccion} applied to \ref{170417_42}. \label{210416_17}
			\item 	$\Gamma, \alpha \inl \beta, \beta \inl \alpha \cl [\beta \inl (\sim \gamma \inl \alpha)] \inl [(\beta \wedge \sim \gamma) \inl \alpha]$  by axiom $(A{\ref{axioma_InfimoAImplicacionVuelta}})$ and  by Lemma \ref{propiedades_Calculo} (\ref{brigAFactores}). \label{210416_18}
			\item 	$\Gamma, \alpha \inl \beta, \beta \inl \alpha \cl (\beta \wedge \sim \gamma) \inl \alpha$  by ($\N$-MP) applied to \ref{210416_17} and \ref{210416_18}. \label{210416_19}
			\item 	$\Gamma, \alpha \inl \beta, \beta \inl \alpha \cl (\beta \wedge \sim \gamma) \inl \sim \gamma$   by axiom $(A{\ref{axioma_infimo_derecha}})$. \label{210416_20}
			\item 	$\Gamma, \alpha \inl \beta, \beta \inl \alpha \cl (\beta \wedge \sim \gamma) \inl (\alpha \wedge \sim \gamma)$  by axiom $(A{\ref{axioma_mayor_cota_inferior}})$ and ($\N$-MP) applied to \ref{210416_19} and \ref{210416_20}. \label{210416_22}
			\item 	$\Gamma, \alpha \inl \beta, \beta \inl \alpha \cl \sim(\beta \il \gamma) \inl (\beta \wedge \sim \gamma)$  by axiom $(A{\ref{axioma_paraSN1}})$. \label{210416_21}
			\item 	$\Gamma, \alpha \inl \beta, \beta \inl \alpha \cl \sim(\beta \il \gamma) \inl (\alpha \wedge \sim \gamma)$  by axiom $(A{\ref{axioma_transitividadNelson}})$ and ($\N$-MP) applied to \ref{210416_21} and \ref{210416_22}. \label{210416_23}
			\item 	$\Gamma, \alpha \inl \beta, \beta \inl \alpha \cl (\alpha \wedge \sim \gamma) \inl \sim (\alpha \il \gamma)$  by axiom $(A{\ref{axioma_paraSN2}})$. \label{210416_24}
			\item 	$\Gamma, \alpha \inl \beta, \beta \inl \alpha \cl \sim(\beta \il \gamma)  \inl \sim (\alpha \il \gamma)$  by axiom $(A{\ref{axioma_transitividadNelson}})$ and ($\N$-MP) applied to \ref{210416_23} and \ref{210416_24}. \label{170417_43}
			\item 	$\Gamma, \alpha \inl \beta \cl (\beta \inl \alpha) \inl [\sim(\beta \il \gamma)  \inl \sim (\alpha \il \gamma)]$  by Theorem \ref{teorema_deduccion} applied to \ref{170417_43}. \label{170417_44}
			\item 	$\Gamma \cl (\alpha \inl \beta) \inl [(\beta \inl \alpha) \inl [\sim(\beta \il \gamma)  \inl \sim (\alpha \il \gamma)]]$  by Theorem \ref{teorema_deduccion} applied to \ref{170417_44}.
		\end{enumerate}
		
		\item[(\ref{120716_03})] 		
		\begin{enumerate}
			\item 	$\Gamma, \sim \alpha \inl \sim \beta,  \sim \beta \inl \sim \alpha \cl \sim \alpha \inl \sim \beta$.
			\item 	$\Gamma, \sim \alpha \inl \sim \beta,  \sim \beta \inl \sim \alpha \cl \gamma \inl (\sim \alpha \inl \sim \beta)$   by Lemma \ref{propiedades_Calculo} (\ref{implicacionATeorema}). \label{220416_01}
			\item 	$\Gamma, \sim \alpha \inl \sim \beta,  \sim \beta \inl \sim \alpha \cl [\gamma \inl (\sim \alpha \inl \sim \beta)] \inl [(\gamma \wedge \sim \alpha) \inl \sim \beta]$  by axiom $(A{\ref{axioma_InfimoAImplicacionVuelta}})$ and  by Lemma \ref{propiedades_Calculo} (\ref{brigAFactores}). \label{220416_02}
			\item 	$\Gamma, \sim \alpha \inl \sim \beta,  \sim \beta \inl \sim \alpha \cl (\gamma \wedge \sim \alpha) \inl \sim \beta$  by ($\N$-MP) applied to \ref{220416_01} and \ref{220416_02}.  \label{220416_04}
			\item 	$\Gamma, \sim \alpha \inl \sim \beta,  \sim \beta \inl \sim \alpha \cl (\gamma \wedge \sim \alpha) \inl \gamma$  by axiom $(A{\ref{axioma_infimo_izquierda}})$.  \label{220416_03}
			\item 	$\Gamma, \sim \alpha \inl \sim \beta,  \sim \beta \inl \sim \alpha \cl (\gamma \wedge \sim \alpha) \inl (\gamma \wedge \sim \beta)$  by axiom $(A{\ref{axioma_mayor_cota_inferior}})$ and ($\N$-MP) applied to \ref{220416_03} and \ref{220416_04}. \label{220416_05}
			\item 	$\Gamma, \sim \alpha \inl \sim \beta,  \sim \beta \inl \sim \alpha \cl \sim(\gamma \il \alpha) \inl (\gamma \wedge \sim \alpha)$  by axiom $(A{\ref{axioma_paraSN1}})$. \label{220416_06}
			\item 	$\Gamma, \sim \alpha \inl \sim \beta,  \sim \beta \inl \sim \alpha \cl \sim(\gamma \il \alpha) \inl  (\gamma \wedge \sim \beta)$  by axiom $(A{\ref{axioma_transitividadNelson}})$ and ($\N$-MP) applied to \ref{220416_05} and \ref{220416_06}. \label{220416_07}
			\item 	$\Gamma, \sim \alpha \inl \sim \beta,  \sim \beta \inl \sim \alpha \cl (\gamma \wedge \sim \beta) \inl \sim(\gamma \il \beta)$  by axiom $(A{\ref{axioma_paraSN2}})$. \label{220416_08}
			\item 	$\Gamma, \sim \alpha \inl \sim \beta,  \sim \beta \inl \sim \alpha \cl  \sim(\gamma \il \alpha) \inl \sim(\gamma \il \beta) $  by axiom $(A{\ref{axioma_transitividadNelson}})$ and ($\N$-MP) applied to \ref{220416_07} and \ref{220416_08}.
			\item 	$\Gamma, \sim \alpha \inl \sim \beta \cl (\sim \beta \inl \sim \alpha) \inl [\sim(\gamma \il \alpha) \inl \sim(\gamma \il \beta)]$  by Theorem \ref{teorema_deduccion}.
			\item 	$\Gamma \cl (\sim \alpha \inl \sim \beta) \inl [(\sim \beta \inl \sim \alpha) \inl [\sim(\gamma \il \alpha) \inl \sim(\gamma \il \beta)]]$   by Theorem \ref{teorema_deduccion}.
		\end{enumerate}
		
		\item[(\ref{condicionIL3d})] 
		\begin{enumerate}[1.]
			\item $\Gamma, \alpha \brig \beta, \beta \brig \alpha, \gamma \brig t, t \brig \gamma \cl  \alpha \brig \beta$.
			\item $\Gamma, \alpha \brig \beta, \beta \brig \alpha, \gamma \brig t, t \brig \gamma \cl  \alpha \inl \beta$ by  \ref{propiedades_Calculo} (\ref{brigAFactores}). \label{111215_46}
			\item $\Gamma, \alpha \brig \beta, \beta \brig \alpha, \gamma \brig t, t \brig \gamma \cl  \sim \beta \inl \sim \alpha$ by \ref{propiedades_Calculo} (\ref{brigAFactores}). \label{211215_17}
			\item $\Gamma, \alpha \brig \beta, \beta \brig \alpha, \gamma \brig t, t \brig \gamma \cl  \beta \brig \alpha$.
			\item $\Gamma, \alpha \brig \beta, \beta \brig \alpha, \gamma \brig t, t \brig \gamma \cl  \beta \inl \alpha$ by  \ref{propiedades_Calculo} (\ref{brigAFactores}). \label{111215_48}
			\item $\Gamma, \alpha \brig \beta, \beta \brig \alpha, \gamma \brig t, t \brig \gamma \cl  \sim \alpha \inl \sim \beta$ by  \ref{propiedades_Calculo} (\ref{brigAFactores}). \label{211215_16}
			\item $\Gamma, \alpha \brig \beta, \beta \brig \alpha, \gamma \brig t, t \brig \gamma \cl  (\alpha \inl \beta) \inl [(\beta \inl \alpha) \inl [(\alpha \il \gamma) \inl (\beta \il \gamma)]]$   by axiom $(A{\ref{axioma_buena_def_implica_der}})$. \label{111215_47}
			\item $\Gamma, \alpha \brig \beta, \beta \brig \alpha, \gamma \brig t, t \brig \gamma \cl  (\beta \inl \alpha) \inl [(\alpha \il \gamma) \inl (\beta \il \gamma)]$  by ($\N$-MP) applied to \ref{111215_46} and \ref{111215_47}. \label{111215_49}
			\item $\Gamma, \alpha \brig \beta, \beta \brig \alpha, \gamma \brig t, t \brig \gamma \cl  (\alpha \il \gamma) \inl (\beta \il \gamma)$  by ($\N$-MP) applied to \ref{111215_48} and \ref{111215_49}. \label{111215_54}
			\item $\Gamma, \alpha \brig \beta, \beta \brig \alpha, \gamma \brig t, t \brig \gamma \cl  (\alpha \inl \beta) \inl [(\beta \inl \alpha) \inl [(\sim (\beta \il \gamma)) \inl (\sim(\alpha \il \gamma))]]$  by 	Lemma \ref{propiedades_Calculo_conDeduccion_para_reticulado} (\ref{120716_02}).
			\label{111215_51}
			\item $\Gamma, \alpha \brig \beta, \beta \brig \alpha, \gamma \brig t, t \brig \gamma \cl  (\beta \inl \alpha) \inl [(\sim (\beta \il \gamma)) \inl (\sim(\alpha \il \gamma))]$  by ($\N$-MP) applied to \ref{111215_46} and \ref{111215_51}. \label{111215_53}
			\item $\Gamma, \alpha \brig \beta, \beta \brig \alpha, \gamma \brig t, t \brig \gamma \cl  (\sim (\beta \il \gamma)) \inl (\sim(\alpha \il \gamma))$  by ($\N$-MP) applied to \ref{111215_48} and \ref{111215_53}.  \label{111215_55}
			\item $\Gamma, \alpha \brig \beta, \beta \brig \alpha, \gamma \brig t, t \brig \gamma \cl  (\alpha \il \gamma) \brig (\beta \il \gamma)$  by  \ref{propiedades_Calculo} (\ref{021115_01})   applied to \ref{111215_54}  and \ref{111215_55}. \label{111215_65}
			\item $\Gamma, \alpha \brig \beta, \beta \brig \alpha, \gamma \brig t, t \brig \gamma \cl  \gamma \brig t$.
			\item $\Gamma, \alpha \brig \beta, \beta \brig \alpha, \gamma \brig t, t \brig \gamma \cl  \gamma \inl t$ by  \ref{propiedades_Calculo} (\ref{brigAFactores}). \label{111215_57}
			\item $\Gamma, \alpha \brig \beta, \beta \brig \alpha, \gamma \brig t, t \brig \gamma \cl  \sim t \inl \sim \gamma$  by  \ref{propiedades_Calculo} (\ref{brigAFactores}). \label{190416_01}
			\item $\Gamma, \alpha \brig \beta, \beta \brig \alpha, \gamma \brig t, t \brig \gamma \cl  t \brig \gamma$.
			\item $\Gamma, \alpha \brig \beta, \beta \brig \alpha, \gamma \brig t, t \brig \gamma \cl  t \inl \gamma$ by  \ref{propiedades_Calculo} (\ref{brigAFactores}).
			\label{111215_58}
			\item $\Gamma, \alpha \brig \beta, \beta \brig \alpha, \gamma \brig t, t \brig \gamma \cl  \sim \gamma \inl \sim t$  by  \ref{propiedades_Calculo} (\ref{brigAFactores}). \label{190416_03}
			\item $\Gamma, \alpha \brig \beta, \beta \brig \alpha, \gamma \brig t, t \brig \gamma \cl  (\gamma \inl t) \inl [(t \inl \gamma) \inl [(\beta \il \gamma) \inl (\beta \il t)]]$   by axiom $(A{\ref{axioma_buena_def_implica_izq}})$. \label{111215_56}
			\item $\Gamma, \alpha \brig \beta, \beta \brig \alpha, \gamma \brig t, t \brig \gamma \cl  (t \inl \gamma) \inl [(\beta \il \gamma) \inl (\beta \il t)]$  by ($\N$-MP) applied to \ref{111215_57} and \ref{111215_56}. \label{111215_59}
			\item $\Gamma, \alpha \brig \beta, \beta \brig \alpha, \gamma \brig t, t \brig \gamma \cl  (\beta \il \gamma) \inl (\beta \il t)$  by ($\N$-MP) applied to \ref{111215_58} and \ref{111215_59}. \label{111215_60}
			\item $\Gamma, \alpha \brig \beta, \beta \brig \alpha, \gamma \brig t, t \brig \gamma \cl  (\sim t \inl \sim \gamma) \inl [(\sim \gamma \inl \sim t) \inl [(\sim (\beta \il t)) \inl (\sim(\beta \il \gamma))]]$  by  \ref{propiedades_Calculo_conDeduccion_para_reticulado} (\ref{120716_03}). \label{190416_02}
			\item $\Gamma, \alpha \brig \beta, \beta \brig \alpha, \gamma \brig t, t \brig \gamma \cl  (\sim \gamma \inl \sim t) \inl [(\sim (\beta \il t)) \inl (\sim(\beta \il \gamma))]$	 by ($\N$-MP) applied to \ref{190416_01} and \ref{190416_02}. \label{190416_04}
			\item $\Gamma, \alpha \brig \beta, \beta \brig \alpha, \gamma \brig t, t \brig \gamma \cl (\sim (\beta \il t)) \inl (\sim(\beta \il \gamma))$  by ($\N$-MP) applied to \ref{190416_03} and \ref{190416_04}. \label{111215_61}	
			\item $\Gamma, \alpha \brig \beta, \beta \brig \alpha, \gamma \brig t, t \brig \gamma \cl  (\beta \il \gamma) \brig (\beta \il t)$  by \ref{propiedades_Calculo} (\ref{021115_01})   applied to \ref{111215_60}  and \ref{111215_61}. \label{111215_66}			
			\item $\Gamma, \alpha \brig \beta, \beta \brig \alpha, \gamma \brig t, t \brig \gamma \cl  (\alpha \il \gamma) \brig (\beta \il t)$  by  \ref{propiedades_Calculo} (\ref{propiedad_Transitividad})  applied to \ref{111215_65} and \ref{111215_66}.
		\end{enumerate}

		\item[(\ref{condicionIL5})] 
		\begin{enumerate}[1.]
			\item $\Gamma, \alpha \cl [(\alpha \wedge \beta) \inl \alpha] \brig [\alpha \inl (\beta \inl \alpha)]$  by axiom $(A{\ref{axioma_InfimoAImplicacion}})$. \label{170417_45}
			\item $\Gamma, \alpha \cl [(\alpha \wedge \beta) \inl \alpha] \inl [\alpha \inl (\beta \inl \alpha)]$  by  Lemma \ref{propiedades_Calculo} (\ref{brigAFactores}) applied to \ref{170417_45}. \label{141215_01}
			\item $\Gamma, \alpha \cl (\alpha \wedge \beta) \inl \alpha$  by axiom $(A{\ref{axioma_infimo_izquierda}})$. \label{141215_02}
			\item $\Gamma, \alpha \cl \alpha \inl (\beta \inl \alpha)$  by ($\N$-MP) applied to \ref{141215_01} and \ref{141215_02}. \label{141215_03}
			\item $\Gamma, \alpha \cl \alpha$. \label{141215_04}
			\item $\Gamma, \alpha \cl \beta \inl \alpha$  by ($\N$-MP) applied to \ref{141215_03} and \ref{141215_04}. \label{141215_09}
			\item $\Gamma, \alpha \cl  [(\alpha \wedge \sim \alpha) \inl \sim \beta] \brig [\alpha \inl (\sim \alpha \inl \sim \beta)]$  by axiom $(A{\ref{axioma_InfimoAImplicacion}})$. \label{170417_46}
			\item $\Gamma, \alpha \cl [(\alpha \wedge \sim \alpha) \inl \sim \beta] \inl [\alpha \inl (\sim \alpha \inl \sim \beta)]$   by   Lemma \ref{propiedades_Calculo} (\ref{brigAFactores}) applied to \ref{170417_46}. \label{141215_05}
			\item $\Gamma, \alpha \cl (\alpha \wedge \sim \alpha) \inl \sim \beta$   by  \ref{propiedades_Calculo_conDeduccion_para_reticulado} (\ref{120716_04}). \label{141215_06}
			\item $\Gamma, \alpha \cl \alpha \inl (\sim \alpha \inl \sim \beta)$  by ($\N$-MP) applied to \ref{141215_05} and \ref{141215_06}. \label{141215_08}
			\item $\Gamma, \alpha \cl \sim \alpha \inl \sim \beta$    by ($\N$-MP) applied to \ref{141215_04} and \ref{141215_08}. \label{141215_10}
			\item $\Gamma, \alpha \cl \beta \brig \alpha$   by  Lemma \ref{propiedades_Calculo} (\ref{021115_01}) applied to \ref{141215_09} and \ref{141215_10}. 
		\end{enumerate}
		
		\item[(\ref{151116_01})] 
		\begin{enumerate}
			\item 	$\Gamma, \alpha \il \beta, \alpha \cl \alpha$.
			\item 	$\Gamma, \alpha \il \beta, \alpha \cl \beta \inl \alpha$ by \ref{propiedades_Calculo} (\ref{implicacionATeorema}). \label{151116_03}
			\item 	$\Gamma, \alpha \il \beta, \alpha  \cl \beta \inl \beta$ by \ref{propiedades_Calculo} (\ref{propiedad_XimplicaX}). \label{151116_04}
			\item 	$\Gamma, \alpha \il \beta, \alpha  \cl \beta \inl (\alpha \wedge \beta)$  by axiom $(A\ref{axioma_mayor_cota_inferior})$ and ($\N$-MP) applied to \ref{151116_03} and \ref{151116_04}. \label{151116_05}
			\item 	$\Gamma, \alpha \il \beta, \alpha  \cl (\alpha \wedge \beta) \inl \beta$ by axiom $(A\ref{axioma_infimo_izquierda})$. \label{151116_06}
			\item 	$\Gamma, \alpha \il \beta, \alpha  \cl (\alpha \il \beta) \inl (\alpha \il (\alpha \wedge \beta))$  by axiom $(A\ref{axioma_buena_def_implica_izq})$ and ($\N$-MP) applied to \ref{151116_05} and \ref{151116_06}. \label{151116_07}
			\item 	$\Gamma, \alpha \il \beta, \alpha  \cl \alpha \il \beta$. \label{151116_08}
			\item 	$\Gamma, \alpha \il \beta, \alpha  \cl \alpha \il (\alpha \wedge \beta)$  by ($\N$-MP) applied to \ref{151116_07} and \ref{151116_08}.
			\item 	$\Gamma, \alpha \il \beta, \alpha  \cl \alpha \inl \beta$. \label{151116_09}
			\item 	$\Gamma, \alpha \il \beta, \alpha   \cl \alpha$. \label{151116_10}
			\item 	$\Gamma, \alpha \il \beta, \alpha \cl \beta$  by ($\N$-MP) applied to \ref{151116_09} and \ref{151116_10}.
			\item 	$\Gamma, \alpha \il \beta \cl \alpha \inl \beta$ by Theorem \ref{teorema_deduccion}.
			\item 	$\Gamma \cl (\alpha \il \beta) \inl (\alpha \inl \beta)$ by Theorem \ref{teorema_deduccion}.
		\end{enumerate}

	\end{itemize}
\end{Proof}


\textbf{Acknowledgements.} We gratefully acknowledge the constructive comments and corrections offered by the referees of \textit{Studia Logica}. This work was partially supported by CONICET (Consejo Nacional de Investigaciones Cient\'ificas y T\'ecnicas, Argentina). 

\bibliographystyle{plain}

\end{document}